\documentclass[aos,nameyear,dvips]{arximspdf}
\usepackage{dcolumn}
\usepackage{graphicx}

\doi{10.1214/10-AOS802}
\volume{38}
\issue{5}
\pubyear{2010}
\firstpage{2723}
\lastpage{2750}

\makeatletter

\newcommand{\gi}{\mathit{gi}}
\newcommand{\NI}{\mathit{NI}}

\newcolumntype{d}[1]{D{.}{.}{#1}}

\newtheorem{lemma}{Lemma}
\newtheorem{theorem}{Theorem}
\newproclaim{remark}{Remark}

\makeatother

\begin{document}
\begin{frontmatter}

\title{Nonparametric estimation of genewise variance for microarray
data\thanksref{T1}}
\runtitle{Genewise variance estimation}

\thankstext{T1}{Supported by NSF Grant DMS-07-14554 and NIH
Grant R01-GM072611.}

\begin{aug}
\author[A]{\fnms{Jianqing} \snm{Fan}\ead[label=e1]{jqfan@princeton.edu}},
\author[B]{\fnms{Yang} \snm{Feng}\ead[label=e2]{fy2158@columbia.edu}} and
\author[C]{\fnms{Yue S.} \snm{Niu}\corref{}\ead[label=e3]{yueniu@math.arizona.edu}}
\runauthor{J. Fan, Y. Feng and Y. S. Niu}
\affiliation{Princeton University, Columbia University and University
of Arizona}
\address[A]{J. Fan\\
Department of Operations Research\\
\quad and Financial Engineering\\
Princeton University\\
Princeton, New Jersey 08544\\
USA\\
\printead{e1}} 
\address[B]{Y. Feng\\
Department of Statistics\\
Columbia University\\
1255 Amsterdam Avenue, 10th Floor\\
New York, New York 10027\\
USA\\
\printead{e2}}
\address[C]{Y. S. Niu\\
Department of Mathematics\\
University of Arizona\\
617 N. Santa Rita Ave.\\
P.O. Box 210089\\
Tucson, Arizona 85721-0089\\
USA\\
\printead{e3}}
\end{aug}

\received{\smonth{9} \syear{2009}}
\revised{\smonth{1} \syear{2010}}

%
\begin{abstract}
Estimation of genewise variance arises from two important applications
in microarray data analysis: selecting significantly differentially
expressed genes and validation tests for normalization of microarray
data. We approach the problem by introducing a two-way nonparametric
model, which is an extension of the famous Neyman--Scott model and is
applicable beyond microarray data. The problem itself poses
interesting challenges because the number of nuisance parameters is
proportional to the sample size and it is not obvious how the variance
function can be estimated when measurements are correlated. In such a
high-dimensional nonparametric problem, we proposed two novel
nonparametric estimators for genewise variance function and
semiparametric estimators for measurement correlation, via solving a
system of nonlinear equations. Their asymptotic normality is
established. The finite sample property is demonstrated by simulation
studies. The estimators also improve the power of the tests for
detecting statistically differentially expressed genes. The methodology
is illustrated by the data from microarray quality control (MAQC)
project.
\end{abstract}

%

\end{frontmatter}

\section{Introduction}\label{sec1}

Microarray experiments are one of widely used technologies nowadays,
allowing scientists to monitor thousands of gene expressions
simultaneously. One of the important scientific endeavors of
microarray data analysis is to detect statistically differentially
expressed genes for downstream analysis [\citet{Cui05},
\citet{Fan04}, \citet{FanRen06},
\citet{StoreyTibs03}, \citet{Tusher01}]. Standard
$t$-test and $F$-test are frequently employed. However, due to the
cost of the experiment, it is common to see a large number of genes
with a small number of replications. Even in customized arrays where
only several hundreds of genes expressions are measured, the number
of replications is usually limited. As a result, we are facing a
high-dimensional statistical problem with a large number of
parameters and a small sample size.

Genewise variance estimation arises at the heart of microarray data
analysis. To select differentially expressed genes among thousands of
genes, the $t$-test is frequently employed with a stringent control
of type I errors. The degree of freedom is usually small due to
limited replications. The power of the test can be significantly
improved if the genewise variance can be estimated accurately. In
such a case, the $t$-test becomes basically a $z$-test. A simple
genewise variance estimator is the sample variance of replicated
data, which is not reliable due to a relatively small number of
replicated genes. They have direct impact on the sensitivity and
specificity of $t$-test [\citet{Cui05}]. Therefore, novel methods
for estimating the genewise variances are needed for improving the
power of the standard $t$-test.

Another important application of genewise variance estimation arises
from testing whether systematic biases have been properly removed
after applying some normalization method, or selecting the most
appropriate normalization technique for a given array. \citet{FanNiu07}
developed such validation tests (see Section~\ref{sec4}), which
require the estimation of genewise variance. The methods of variance
estimation, like pooled variance estimator, and REML estimator
[\citet{Smyth05}], are not accurate enough due to
the small number of replications.

Due to the importance of genewise variance in microarray data
analysis, conscientious efforts have been made to accurately
estimate it. Various methods have been proposed under different
models and assumptions. It has been widely observed that genewise
variance is to a great extent related to the intensity level.
\citet{Kamb} proposed a crude regression estimation of variance from
microarray control data. \citet{Tong07} discussed a family of
shrinkage estimators to improve the accuracy.

Let $R_{\gi}$ and $G_{\gi}$, respectively, be the intensities of red (Cy3)
and green (Cy5) channels for
the $i$th replication of the $g$th gene on a two-color microarray
data. The log-ratios and
log-intensities are computed, respectively, as
\begin{eqnarray}
Y_{\gi} = \log_2 (G_{\gi}/R_{\gi})\quad \mbox{and}\quad X_{\gi} =
\tfrac{1}{2}\log_2 (G_{\gi}R_{\gi}),\nonumber\\
\eqntext{i = 1,\ldots, I, g=1,\ldots, N,}
\end{eqnarray}
where $I$ is the number of replications for each gene and $N$ is the
number of genes with replications. For the purpose of estimating
genewise variance, we assume that there is no systematic biases or
the systematic biases have been removed by a certain normalization
method. This assumption is always made for selecting significantly
differentially
expressed
genes or validation test under the null hypothesis. Thus, we have
\[
Y_{\gi}=\alpha_g+\sigma_{\gi} \epsilon_{\gi}
\]
with $\alpha_g$ denoting the log-ratio of gene expressions in the
treatment and control samples. Here, $(\epsilon_{g1},\ldots,
\epsilon_{gI})^T$ follows a multivariate normal distribution with
$\epsilon_{\gi} \sim N(0,1)$ and
$\operatorname{Corr}(\epsilon_{\gi},\epsilon_{gj})=\rho$ when $i\neq j$. It
is also assumed that observations from different genes are
independent.
Such a model was used in \citet{Smyth05}.

In the papers
by \citet{Wang08} and \citet{Carroll08}, nonparametric
measurement-error models have been introduced to aggregate the
information of estimating the genewise variance:
%
%
\begin{eqnarray}\label{a0}
Y_{\gi} &=& \alpha_g + \sigma(\alpha_g)
\varepsilon_{\gi},\nonumber\\[-8pt]\\[-8pt]
\operatorname{corr}(\varepsilon_{\gi},\varepsilon_{\gi'}) &=& 0,\qquad g=1,\ldots, N,
i=1,\ldots, I.\nonumber
\end{eqnarray}
The model is intended for the analysis of the Affymetrix array
({one-color array}) data in which
$\alpha_g$ represents the expected intensity level, {and $Y_{\gi}$ is
the $i$th replicate of
observed expression level of gene $g$}. When it is applied to the
two-color microarray data as in our setting, in which $\alpha_g$ is
the relative expression profiles between the
treatment and control, several drawbacks emerge: (a) the model is
difficult to interpret as the
genewide variance is a function of the log-ratio of expression
profiles; (b) errors-in-variable
methods have a very slow rate of convergence for the nonparametric
problem and the observed
intensity information $X_{\gi}$ is not used; (c)~they are usually hard
to be implemented robustly
and depend sensitively on the distribution of $\sigma(\alpha_g)
\varepsilon_{\gi}$ and the
i.i.d. assumption on the noise; 
(d) in many microarray applications, $\alpha_g=0$ for most $g$ and
hence $\sigma(\alpha_g)$ are
the same for most genes, which is unrealistic. Therefore, our model
(\ref{a1}) below is
complementary to that of \citet{Wang08} and \citet
{Carroll08}, with focus on the applications to two-color
microarray data.

To overcome these drawbacks in the applications to microarray data and
to utilize the observed
intensity information, we assume that $\sigma_{\gi} = \sigma(X_{\gi})$
for a smooth function
$\sigma(\cdot)$. This leads to the following two-way nonparametric
model:
%
%
\begin{equation}\label{a1}
Y_{\gi}=\alpha_g+\sigma(X_{\gi})\epsilon_{\gi},\qquad g=1,\ldots, N,
i=1,\ldots, I,
\end{equation}
for estimating genewise variance. This model is clearly an
extension of the Neyman--Scott problem [\citet{Neyman}], in
which the genewise variance is a constant. The Neyman--Scott problem
has many applications in astronomy. Note that the number of
nuisance parameters $\{\alpha_g\}$ is proportional to the sample
size. This imposes an important challenge to the nonparametric
problem. It is not even clear whether the function $\sigma(\cdot)$
can be consistently estimated.

To estimate the genewise variance in their microarray data analysis,
\citet{Fan04} assumed a model similar to (\ref{a1}). But in the
absence of other available techniques, they had to impose that the
treatment effect $\{\alpha_g\}$ is also a smooth function of the
intensity level so that they can apply nonparametric methods to
estimate genewise variance [\citet{Ruppert97}]. However, this
assumption is not valid in most microarray applications, and the
estimator of genewise variance incurs big biases unless
$\{\alpha_g\}$ is sparse, a situation that \citet{Fan04} hoped.
\citet{FanNiu07} approached this problem in another\vspace*{1pt} simple way. When
the noise in the replications is small, that is, $X_{\gi} \approx\bar{X}_g$,
where $\bar{X}_g$ is the sample mean for the $g$th gene. Therefore,
they simply smoothed the pair $\{(\bar{X}_g, \bar{r}_g)\}$, where
$\bar{r}_g = \sum_{i=1}^I (Y_{\gi}-\bar{Y}_g)^2/(I-1)$. This also
leads to a biased estimator, which is denoted as $\hat\xi^2(x)$. One asks
naturally whether the function $\sigma(\cdot)$ is estimable and
how it can be estimated in the general two-way nonparametric model.

We propose a novel nonparametric approach to estimate the genewise
variance. We first study a benchmark case when there is no
correlation between replications, that is, $\rho= 0$. This corresponds
to the case with independent replications across arrays
[\citet{Fan05}, \citet{Huang05}]. It is
also applicable to those dealt by the Neyman--Scott problem. By
noticing $\mathrm{E}\{(Y_{\gi}-\bar{Y}_g)^2 | X_{\gi}\}$ is a linear
combination of $\sigma^2(X_{\gi})$, we obtain a system of linear
equations. Hence, $\sigma^2(\cdot)$ can be estimated via
nonparametric regression of a proper linear combination of
$\{(Y_{\gi}-\bar{Y}_g)^2, i = 1,\ldots, I\}$ on $\{X_{\gi}\}$. The
asymptotic normality of the estimator is established. In the case
that the replication correlation does not vanish, the system of
equations becomes nonlinear and cannot be analytically solved. However,
we are able to derive the correlation corrected estimator, based on
the estimator without genewise correlation. The genewise variance
function and the correlation coefficient of repeated measurements
are simultaneously estimated by iteratively solving a nonlinear
equation. The asymptotic normality of such estimators is
established.

Model (\ref{a1}) can be applied to the microarrays in which
within-array replications are not available. In that case, we can
aggregate all the microarrays together and view them as a super
array with replications [\citet{Fan05}, \citet{Huang05}].
In other words,
$i$ in (\ref{a1}) indexes arrays and $\rho$ can be taken as 0, namely
(\ref{a1}) is the across-array replication with $\rho= 0$.

The structure of this paper is as follows. In Section \ref{sec2}, we discuss
the estimation schemes of the genewise variance and establish the
asymptotic properties of the estimators. Simulation studies are
given in Section \ref{sec3} to verify the performance of our methods in the
finite sample. Applications to the data from Microarray Quality
Control (MAQC) project are showed in Section \ref{sec4} to illustrate the
proposed methodology. In Section \ref{sec5}, we give a short summary. Technical
proofs are relegated to the
\hyperref[app]{Appendix}.

\section{Nonparametric estimators of genewise variance}\label{sec2}

\subsection{Estimation without correlation}\label{sec21}

We first consider the specific case where there is no correlation
among the replications $Y_{g1},\ldots, Y_{gI}$ of
the same gene $g$ under model (\ref{a1}). This is usually applicable
to the across-array replication and stimulates our
procedure for the more general case with the replication
correlation. In the former case, we have
\[
\mathrm{E}[(Y_{\gi}-\bar{Y_g})^2 |{\mathbf{X}}] =
(I-1)^2\sigma^2(X_{\gi})/I^2+\sum_{j \neq i}
\sigma^2(X_{gj})\big/I^2,\qquad i = 1,\ldots, I.
\]
We will discuss in Section \ref{sec224} the case that $I=2$. For $I > 2$, we have
$I$ different equations with $I$ unknowns $\sigma^2(X_{g1})$,
$\sigma^2(X_{g2}), \ldots, \sigma^2(X_{gI})$ for a given gene
$g$. Solving these $I$ equations, we can express the unknowns in
terms of $\{\mathrm{E}[(Y_{\gi}-\bar{Y}_g)^2|{\mathbf{X}}]\}
_{i=1}^I$, estimable
quantities. Let
\[
{\mathbf{r}_g} = \bigl((Y_{g1}-\bar{Y}_g)^2,\ldots, (Y_{gI}-\bar{Y}_g)^2\bigr)^T
\quad\mbox{and}\quad \bolds{\sigma}^2_g = (\sigma^2(X_{g1}),\ldots,
\sigma^2(X_{gI}))^T.
\]
Then, it can easily be shown that
$\bolds{\sigma}^2_g = \mathbf{B} \mathrm{E}[{\mathbf
{r}_g}|{\mathbf{X}}]$, where
$\mathbf{B}$ is the
coefficient matrix:
\[
\mathbf{B}= \bigl((I^2-I){\mathbf{I}} - {\mathbf{E}}\bigr)/(I-1)(I-2)
\]
with ${\mathbf{I}}$ being the $I
\times I$ identity matrix and ${\mathbf{E}}$ the $I \times I$ matrix with
all elements 1. Define
\[
{\mathbf{Z}_g} = (Z_{g1},\ldots, Z_{gI})^T \stackrel{\triangle}{=}
\mathbf{Br}_{\mathbf{g}}.
\]
Then we have
%
%
\begin{equation} \label{b1}
\sigma^2(X_{\gi})=\mathrm{E}[Z_{\gi}|{\mathbf{X}}].
\end{equation}
Note that the left-hand side of (\ref{b1}) depends only on $X_{\gi}$,
not other variables. By the the double expectation formula, it
follows that the variance function $\sigma^2(\cdot)$ can be expressed
as the univariate regression
%
%
\begin{equation}\label{b2}
\sigma^2(x)=\mathrm{E}[Z_{\gi}|X_{\gi}=x],\qquad i = 1,\ldots, I.
\end{equation}
Using the synthetic data $\{(X_{\gi}, Z_{\gi}), g=1,\ldots, N\}$ for
each given $i$, we can apply the local linear regression technique
[\citet{Fan96}] to obtain a nonparametric estimator $\hat{\eta}_i^2
(x)$ of $\sigma^2(\cdot)$. Explicitly, for a given kernel $K$ and
bandwidth~$h$,
%
%
\begin{equation}\label{b3}
\hat{\eta}_i^2 (x) = \sum_{g=1}^N
W_{N,i} \biggl(\frac{X_{\gi}-x}{h} \biggr)Z_{\gi},\qquad i=1,\ldots, I,
\end{equation}
with
\[
W_{N,i}(u) = h^{-1} K(u)\frac{S_{N,2}-u S_{N,1}}{S_{N,2} S_{N, 0} - S_{N,
1}^2},
\]
where $K_h(u) = h^{-1} K(u/h)$ and $S_{N,l} = \sum_{g=1}^N
K_h(X_{\gi} - x)[(X_{\gi}-x)/h ]^l$, whose
dependence on $i$ is suppressed. Thus,\vspace*{1pt} we have $I$ estimators
$\hat{\eta}^2_1(x), \ldots, \hat{\eta}^2_I(x)$ for the same
genewise variance function $\sigma^2(\cdot)$. Each of these $I$
estimators $\hat\eta^2_i(x)$ is a consistent estimator of
$\sigma^2(x)$. To optimally\vspace*{2pt} aggregate those $I$ estimators, we need
the asymptotic properties of $\bolds{\eta}(x) = (\hat{\eta}^2_1(x), \ldots, \hat{\eta}^2_I(x))^T$.

Denote
\begin{eqnarray*}
c_K&=&\int_{-\infty}^{\infty}u^2K(u)\,du,\qquad
d_K=\int_{-\infty}^{\infty}K^2(u)\,du,\\
\sigma_1&=&\mathrm{E}[\sigma(X_{\gi})]\quad \mbox{and}\quad \sigma_2=\mathrm
{E}[\sigma
^2(X_{\gi})].
\end{eqnarray*}
Assume that $X_{\gi}$ are i.i.d. with marginal density
$f_X(\cdot)$ and $\varepsilon_{\gi}$ are i.i.d. random variables from
the standard normal distribution. In the following result, we assume
that $I$ is fixed, but $N$ diverges.
\begin{theorem}\label{T1}
Under the regularity conditions in the \hyperref[app]{Appendix}, for a fixed point
$x$, we have
\[
\bolds\Sigma^{-1/2}\bigl(\bolds{\eta}-\bigl(\sigma^2(x)+b(x)+o_P(h^2)\bigr)
\mathbf{e}\bigr)
\stackrel{D}\longrightarrow N(\mathbf{0}, \mathbf{I}),
\]
provided that $h \rightarrow0$ and $Nh \rightarrow\infty$, where
$\mathbf{e}= (1, 1 ,\ldots, 1)^T$ and
\[
\bolds\Sigma= V_1 {\mathbf{I}} + V_2 ({\mathbf{E}- \mathbf{I}})
\]
with $ b(x) = \frac{h^2}{2}c_K(\sigma^2(x))''$,
\begin{eqnarray*}
V_1 &=& \frac{d_K}{Nhf_X(x)} \biggl\{2\sigma^4(x)+ \frac
{4+4(I-1)(I-3)}{(I-1)(I-2)^2}\sigma_2\sigma^2(x)+\frac
{2}{(I-1)(I-2)}\sigma_2^2 \biggr\},\\
V_2 &=& \frac{1}{N} \biggl\{\frac{4}{(I-1)^2}\sigma^4(x)-\frac
{8}{(I-1)^2}\sigma_2\sigma^2(x) + \frac{2(I-3)}{(I-1)^2(I-2)}\sigma
^2_2 \biggr\}.
\end{eqnarray*}
\end{theorem}

Note that $V_2$ is one order of magnitude smaller than $V_1$. Hence,
the estimators $\hat\eta^2_1(x),\ldots, \hat\eta^2_I(x)$ are
asymptotically independently distributed as $ N(\sigma^2(x) + b(x),
V_1)$. Their dependence is only in the second order. The best
linear combination of $I$ estimators is
%
%
\begin{equation} \label{b4}
\hat\eta^2(x)=[\hat{\eta}^2_1(x)+ \hat{\eta}^2_2(x)+ \cdots+
\hat{\eta}^2_I(x)]/I
\end{equation}
with the asymptotic distribution
%
%
\begin{equation} \label{b5}
N \bigl( \sigma^2(x) + b(x), V_1/I + (1-1/I)V_2 \bigr).
\end{equation}
See also the aggregated estimator (\ref{b12}) with $\rho= 0$, which
has the same asymptotic property as the estimator (\ref{b6}). See
Remark \ref{remark1} below for additional discussion.

Theorem \ref{T1} gives the asymptotic normality of the proposed
nonparametric estimators under the
presence of a large number of nuisance parameters $\{\alpha_g\}
_{g=1}^N$. With the newly proposed
technique, we do not have to impose any assumptions on $\alpha_g$ such
as sparsity or smoothness.
This kind of local linear estimator can be applied to most two-color
microarray data, for instance,
customized arrays and Agilent arrays.



\subsection{Variance estimation with correlated replications}\label{sec22}

\subsubsection{Aggregated estimator}

We now consider the case with correlated with-array replications.
There is a lot of evidence that correlation among within-array
replicated genes exists [\citet{Smyth05}, \citet
{FanNiu07}]. Suppose that
within-array replications have a common correlation
$\operatorname{corr}(Y_{\gi},Y_{gj}|{\mathbf{X}}) = \rho$ when $i
\neq j$. Observations
across different genes or arrays are independent. Then the
conditional variance of $(Y_{\gi}-\bar{Y}_g)$ can be expressed as
%
%
\begin{eqnarray}\label{b6}
&&\operatorname{var}[(Y_{\gi}-\bar{Y_g})|{\mathbf{X}}] \nonumber\\
&&\qquad= (I-1)^2\sigma
^2(X_{\gi})/I^2 +
2\rho\mathop{\sum_{1\leq j< k \leq I,}}_{j\neq i, k\neq i}\sigma
(X_{gj})\sigma(X_{gk})\big/I^2\\
&&\qquad\quad{} + 2(I-1)\rho\sum_{j\neq i}\sigma^2(X_{gj})\big/I^2-\sum_{j\neq
i}\sigma(X_{\gi})\sigma(X_{gj})\big/I^2.\nonumber
\end{eqnarray}
This is a complex system of nonlinear equations and the analytic
form cannot be found. Innovative ideas are needed.

Using the same notation as that in the previous section, it can be
calculated that
\begin{eqnarray*}
\mathrm{E}[Z_{\gi}|{\mathbf{X}}] &=& \sigma^2(X_{\gi})-\frac
{2}{I-1}\sum_{j\neq i}
\rho\sigma(X_{\gi})\sigma(X_{gj})\\
&&{} +\frac{2}{(I-1)(I-2)}\mathop{\sum_{1\leq j< k \leq I,}}_{j\neq i, k\neq
i} \rho\sigma(X_{gj})\sigma(X_{gk}).
\end{eqnarray*}
Taking the expectation with respect to $X_{gj}$ for all $j \not=
i$, we obtain
%
%
\begin{equation}\label{b7}
\mathrm{E}[Z_{\gi}|{X_{\gi}=x}] = \sigma^2(x) - 2\rho\sigma_1\sigma
(x) +
\rho\sigma_1^2 \stackrel{\triangle}{=} \eta^2(x),
\end{equation}
where $\sigma_1 = \mathrm{E}[\sigma(X)]$.


Here, we can directly apply the local linear approach to all
aggregated data $\{(X_{\gi},Z_{\gi})\}_{i,g=1}^{I,N}$, due to the same
regression function (\ref{b7}). Let $\hat\eta_A^2(\cdot)$ be the
local linear estimator of $\eta^2(\cdot)$, based on the aggregated
data. Then
%
%
\begin{equation}\label{b15}
\hat\eta_A^2(x) = \sum_{g=1}^N \sum_{i=1}^I W_{N}
\biggl(\frac{X_{\gi}-x}{h} \biggr) Z_{\gi}
\end{equation}
with
\[
W_{N} (u) = h^{-1} K(u)\frac{S_{\NI,2}-u S_{\NI,1}} {S_{\NI, 0}
S_{\NI,2}-S_{\NI,1}^2},
\]
where
$
S_{\NI,l} = \sum_{g=1}^N \sum_{i=1}^I K_h ( X_{\gi}-x )
[ ( X_{\gi}-x)/h]^l.
$
There are two solutions to~(\ref{b7}):
%
%
\begin{equation}\label{b19}
\hat{\sigma}_A(x,\rho)^{(1),(2)} = \hat\rho\hat\sigma_1 \pm
\sqrt{\hat\rho^2\hat\sigma_1^2 -
\hat\rho\hat\sigma_1^2+\hat\eta_A^2(x)},
\end{equation}
Notice that given the sample ${\mathbf{X}}$ and ${\mathbf{Y}}$,
$\hat\sigma_A(x,\rho)^{(1),(2)}$ are continuous in both $x$
and~$\rho$. For $\rho<0$, $\hat\sigma_A(x,\rho)^{(1)}$ should be used
since the standard deviation should be nonnegative. Since
$\hat\sigma_A(x,\rho)^{(1)}>\hat\sigma_A(x,\rho)^{(2)}$ for every
$x$ and $\rho$, by the continuity of the solution in $\rho$, we can
only use the same solution when $\rho$ changes continuously. Then
$\hat\sigma_A(x,\rho)^{(1)}$ should always be used regardless of
$\rho$. From now on, we drop the superscript and denote
%
%
\begin{equation}\label{b8}
\hat{\sigma}_A(x) = \rho\sigma_1+\sqrt{\rho^2\sigma_1^2-\rho
\sigma_1^2+\hat\eta^2_A(x)}.
\end{equation}
This is called the aggregated estimator.
Note that in (\ref{b8}), $\rho$, $\sigma_1$ and $\sigma(\cdot)$
are all
unknown.

\subsubsection{Estimation of correlation}\label{sec222}
To estimate $\rho$, we
assume that there are $J$ independent arrays ($J\geq2$). In other
words, we
observed data from (\ref{a1}) independently $J$ times. In this
case, the residual maximum likelihood (REML) estimator introduced by
\citet{Smyth05} is as follows:
%
%
\begin{equation} \label{b9}
\hat{\rho}_0= \frac{\sum_{g=1}^N s_{B, g}^2 - \sum_{g=1}^N s_{W, g}^2
}{\sum_{g=1}^N s_{B, g}^2 + (I-1)\sum_{g=1}^N s_{W, g}^2} ,
\end{equation}
where $s_{B, g}^2= I(J-1)^{-1} \sum_{j=1}^J
(\bar{Y}_{gj}-\bar{Y}_{g})^2$ with $\bar{Y}_{gj} = I^{-1}
\sum_{i=1}^I Y_{\mathit{gij}}$ and $\bar Y_g =J^{-1} \sum_{j=1}^J
\bar{Y}_{gj}$ is the between-arrays variance and $s_{W,g}^2$ is the
within-array variance:
\[
s_{W,g}^2=\frac{1}{J(I-1)}\sum_{j=1}^J\sum_{i=1}^I(Y_{\mathit{gij}}-\bar{Y}_{gj})^2.
\]
As discussed in \citet{Smyth05}, the estimator $\hat\rho_0$ of
$\rho$ is consistent when $\operatorname{var}(Y_{\mathit{gij}} | \mathbf{X})
= \sigma_{g}$ is
the same for all $i=1,\ldots,I$ and $j=1,\ldots,J$. However, this
assumption is not valid under the model (\ref{a1}) and a correction
is needed. We propose the following estimator:
%
%
\begin{equation} \label{b10}
\hat\rho= \frac{\sigma_2}{\sigma_1^2} \cdot\frac{\sum_{g=1}^N
s_{B, g}^2 - \sum_{g=1}^N s_{W, g}^2
}{\sum_{g=1}^N s_{B, g}^2 + (I-1)\sum_{g=1}^N s_{W, g}^2}.
\end{equation}
The consistency of $\hat\rho$ is given by the following theorem.
\begin{theorem}\label{P1}
Under the regularity condition in the \hyperref[app]{Appendix}, the estimator $\hat
\rho$ of $\rho$ is $\sqrt{N}$-consistent:
\[
\hat\rho- \rho= O_P(N^{-1/2}).
\]
\end{theorem}

With a consistent estimator of $\rho$, $\sigma_1$, $\sigma_2$ and
$\sigma_A(\cdot)$ can be solved by the following iterative algorithm:
\begin{enumerate}
\item[Step 1.] Set $\hat\eta_A^2(\cdot)$ as an initial estimate of
$\sigma^2_A(\cdot)$.
\item[Step 2.] With $\hat\sigma_A(\cdot)$, compute
%
%
\begin{equation} \label{b11}\quad
\hat{\sigma}_1 = N^{-1} \sum_{g=1}^N \hat{\sigma}_A(X_{\gi}),\qquad
\hat{\sigma}_2 = N^{-1} \sum_{g=1}^N \hat{\sigma}^2_A(X_{\gi}),\qquad
\hat\rho=\hat\rho_0\hat{\sigma}_2/\hat{\sigma}_1^2 .
\end{equation}
\item[Step 3.] With $\hat{\sigma}_1$, $\hat{\sigma}_2$ and $\hat
\rho$, compute $\hat\sigma_A(\cdot)$ using (\ref{b8}).
\item[Step 4.] Repeat steps 2 and 3 until convergence.
\end{enumerate}
This provides simultaneously the estimators $\hat{\sigma}_1$, $\hat
{\sigma}_2$, $\hat\rho$ and $\hat{\sigma}_A(\cdot)$.
From our numerical experience, this algorithm converges quickly after a
few iterations.
When the algorithm converges, the estimator $\sigma^2_A(x)$ is given by
%
%
\begin{equation}\label{b12}
\hat\sigma_A(x) =\hat\rho\hat\sigma_1+\sqrt{\hat\rho^2
\hat\sigma_1^2-\hat\rho\hat\sigma_1^2+\hat\eta^2_A(x)} .
\end{equation}

Note that the presence of multiple arrays is only used to estimate the
correlation $\rho$ for the replications. It is not needed for
estimating the genewise variance function. In the case of the
presence of $J$ arrays, we can take the average of the
$J$ estimates from each array.

\subsubsection{Asymptotic properties}

Following a similar idea as the case without correlation, we can
derive the asymptotic property of $\hat\eta^2_A(x)$.
\begin{theorem}\label{T3}
Under the regularity conditions in the \hyperref[app]{Appendix}, for a fixed point $x$,
we have
\[
\{V^{*}\}^{-1/2} \{\hat\eta_A^2(x) - [\eta^2(x) + \beta
(x)]+o_P(h^2) \}
\stackrel{D} \longrightarrow N(0,1),
\]
provided that $h \rightarrow0$ and $Nh \rightarrow\infty$, with
$\beta(x) = \frac{h^2}{2}c_K(\eta^2(x))''$ and
\[
V^{*} = \frac{1}{I}V'_1+\frac{I-1}{I}V'_2,
\]
where
\begin{eqnarray*}
V'_1 &=& \frac{d_K}{Nhf_X(x)} \{ 2\sigma^4(x)-8\rho\sigma_1\sigma
^3(x) +C_2 \sigma^2(x) + C_3 \sigma(x) + C_4 \},\\
V'_2 &=& \frac{1}{N} \{ D_0\sigma^4(x) + D_1\sigma^3(x) +
D_2\sigma^2(x) +D_3\sigma(x) + D_4 \}
\end{eqnarray*}
with coefficients $C_2,\ldots, C_4, D_0,\ldots, D_4$ defined in
the \hyperref[app]{Appendix}.
\end{theorem}

The asymptotic normality of $\hat\sigma_A^2(x)$ can be derived from
that of $\hat\eta_A^2(x)$. More specifically, $\hat\sigma^2_A(x) =
\varphi(\eta_A^2(x))$ with $\varphi(z) = (\rho\sigma_1 +
\sqrt{\rho^2\sigma_1^2-\rho\sigma_1^2+z} )^2$. The derivative of
$\varphi(\cdot)$ with respect to $z$ is $\psi(z) = \rho\sigma_1 /
\sqrt{\rho^2\sigma^2_1-\rho\sigma_1^2+z} + 1$. Then, by the delta
method, we have
\[
\{V^{*}\}^{-1/2} \bigl(\hat\sigma_A^2(x) - \varphi\bigl(\eta^2(x)
+ \beta(x)+o_P(h^2)\bigr)\bigr)\stackrel{D}\longrightarrow
N(0,\psi^2(\eta^2(x))).
\]
\begin{remark}\label{remark1}
An alternative approach when correlation exists is to apply the same
correlation correction idea to $\{X_{\gi},Z_{\gi}\}^N_{g=1}$ for every
replication $i$, resulting in the estimator $\hat{\sigma}_i^2(x)$. In
this case, it can be proved
that the best linear combination of the estimator is
%
%
\begin{equation} \label{b13}
\hat\sigma^2(x)=[\hat{\sigma}^2_1(x)+\hat{\sigma}^2_2(x)+\cdots+
\hat{\sigma}^2_I(x)]/I.
\end{equation}
This estimator has the same asymptotic performance as the aggregated
estimator. However, we prefer the aggregated estimator due to the
following reasons: the equation (\ref{b12}) only needs to be solved
once by using the algorithm in Section \ref{sec222}, all data are treated
symmetrically, and $\hat{\eta}_A^2(\cdot)$ can be estimated more
stably.
\end{remark}

\subsubsection{Two replications}\label{sec224}

The aforementioned methods apply to the case when there are more than two
replications. For the case $I=2$, the equations for
$\operatorname{var}[(Y_{\gi}-\bar{Y}_g)|{\mathbf{X}}]$ collapse into
one. In this case,
it can be shown using the same arguments before that
%
%
\begin{equation} \label{b17}\quad
\operatorname{var}[(Y_{\gi}-\bar{Y}_g)|X_{\gi} =x] =
\tfrac{1}{4}\sigma^2(x) + \tfrac{1}{4}\sigma_2
-\tfrac{1}{2}\rho\sigma_1\sigma(x),\qquad i = 1, 2,
\end{equation}
where $\sigma_2 = \mathrm{E}[\sigma^2(X_{\gi})]$. In this case, the
left-hand side is always equal to $\operatorname{var}[(Y_{g1} -
Y_{g2})/2|X_{\gi} =x]$.

Let $\hat\eta^2(x)$ be the local linear estimator of the function on
the right-hand side by smoothing $\{(Y_{g1}-Y_{g2})^2/4\}_{g=1}^N$
on $\{X_{g1}\}_{g=1}^N$ and $\{X_{g2}\}_{g=1}^N$. Then the genewise
variance is a solution to the following equation:
%
%
\begin{equation} \label{b18}
\hat\sigma(x) = \hat\rho\hat\sigma_1 + \sqrt{\hat\rho^2\hat
\sigma^2_1
-\hat\sigma_2 + 4\hat\eta^2(x)}.
\end{equation}
The algorithm in Section \ref{sec222} can be applied directly.


\section{Simulations and comparisons}\label{sec3}

In this section, we conduct simulations to evaluate the finite
sample performance of different variance estimators $\hat\xi^2(x)$,
$\hat\eta^2(x)$
and $\hat\sigma_A^2(x)$. First, the bias problem of the naive
nonparametric variance estimator $\hat\xi^2(x)$ is demonstrated. It is shown
that this bias
issue can be eliminated by our newly proposed methods. Then we consider
the estimators $\hat\eta^2(x)$ and $\hat\sigma_A^2(x)$ under different
configurations of the within-array replication correlation.

\subsection{Simulation design}\label{sec31}

In all the simulation examples, we set the number of genes $N =
2000$, each gene having $I = 3$ within-array replications and $J=4$
independent arrays. For the purpose of investigating the genewise
variance estimation, the data are generated from model (\ref{a1}).
The details of simulation scheme are summarized as follows:
\begin{itemize}[$\alpha_g$:]
\item[$\alpha_g$:] The expression levels of the first 250 genes are
generated from the standard double exponential distribution. The
rest are 0s. These expression levels are the same over 4 arrays
in each simulation, but may vary over simulations.
\item[$X$:] The intensity is generated from a mixture
distribution: with probability 0.7 from the distribution $0.0004
(x-6)^3 I (6<x<16)$ and 0.3 from the uniform distribution
over $[6, 16]$.
\end{itemize}

\begin{itemize}[$\sigma^2(\cdot)$:]
\item[$\varepsilon$:] $\varepsilon_{\gi}$ is generated from the
standard normal
distribution.
\item[$\sigma^2(\cdot)$:] The genewise variance function is taken as
\[
\sigma^2(x)=0.15 +0.015 (12-x)^2 I\{x<12\}.
\]
\end{itemize}
The parameters are taken from \citet{Fan05}. The kernel function is
selected as $\frac{70}{81}(1-|x|^3)^3I(|x| \leq1)$. In addition, we
fix the
bandwidth $h=1$ for all the numerical analysis.

For every setting, we repeat the whole simulation process for
$T$ times and evaluate the estimates of $\sigma^2(\cdot)$ over
$K=101$ grid points $\{x_k\}_{k=1}^K$ on the interval $[6,16]$.
For the $k$th grid point, we define
\begin{eqnarray*}
B_k&=&\bar{\sigma}^2(x_k) - \sigma^2(x_k) \qquad\mbox{with }
\bar{\sigma}^2(x_k) = T^{-1}\sum_{t=1}^T\hat\sigma_t^2(x_k),
\\
S_k&=&T^{-1}\sum_{t=1}^T[\hat\sigma_t^2(x_k)-\bar{\sigma}^2(x_k)]^2,
\end{eqnarray*}
%
and $\mathrm{MSE}_k=B_k^2+S_k$.
Let
$f(\cdot)$ be the density function of intensity $X$. Let
\[
\mathrm{Bias}^2 = \sum_{k=1}^KB_k^2f(x_k) \Big/
\sum_{k=1}^Kf(x_k),\qquad
\mathrm{VAR}= \sum_{k=1}^KS_kf(x_k) \Big/ \sum_{k=1}^Kf(x_k)
\]
and
\[
\mathrm{MISE}=\sum_{k=1}^K\mathrm{MSE}_kf(x_k) \Big/\sum_{k=1}^Kf(x_k)
\]
be the integrated squared bias ($\mathrm{Bias}^2$), the integrated
variance (VAR), and
the integrated mean squared error (MISE) of the estimate
$\hat\sigma^2(\cdot)$, respectively. For the $t$th simulation
experiment, we define
\[
\mathrm{ISE}_t=\sum_{k=1}^{K}\bigl(\hat\sigma_t^2(x_k)-\sigma^2(x_k)\bigr)^2f(x_k)
\Big/ \sum_{k=1}^Kf(x_k)
\]
be the integrated squared error for the $t$th simulation.

\subsection{The bias of naive nonparametric estimator}\label{sec32}

A naive approach is to regard $\alpha_g$ in (\ref{a1}) as a smooth
function of $X_{\gi}$, namely, $\alpha_g = \alpha(X_{\gi})$. The
function $\alpha(\cdot)$ can be estimated by a local linear
regression estimator, resulting in an estimated function
$\hat{\alpha}(\cdot)$. The squared residuals $\{r_{\gi}^2\}_{g=1}^N$
is then further smoothed on $\{X_{\gi}\}_{g=1}^N$ to obtain an
estimate $\hat{\xi}^2(x)$ of the variance function
${\sigma}^2(\cdot)$, where $ r_{\gi} = \hat{Y}_{\gi} -
\hat{\alpha}(X_{\gi})$ [\citet{Ruppert97}].

To provide a comprehensive view of the performances of the
naive and the new estimators, we first compare the performances of
$\hat\xi^2(x)$ and
$\hat\eta^2(x)$ under the smoothness assumption of the gene effect
$\alpha
_g$. Data from the naive nonparametric
regression model is also generated with
\[
\alpha(x) = \exp\biggl( -\frac{1}{1-(x-13)^2} \biggr)I\{12< x <
14\}.
\]
This allows us to understand the loss of efficiency when $\alpha_g$ is
continuous in $X_{\gi}$. This usually does not occur for microarray
data, but can appear in other applications.
Note that $\alpha(\cdot)$ is zero in most of the region and thus is
reasonably sparse. Here, the number of simulations is taken to be
$T=100$. The data is generated with the assumption that $\rho= 0$,
in which case the variance estimators $\hat\eta^2(x)$ and $\hat\sigma
_A^2(x)$ have
the same performance (see also Table \ref{Tb:Comp_NC_WC_AGG} below). Thus, we only report
the performance of $\hat\eta^2(x)$.

%
\begin{table}
\caption{Mean integrated squared bias ($\mathrm{Bias}^2$), mean
integrated variance
(VAR), mean integrated squared error (MISE) over 100 simulations for
variance estimators $\hat\xi^2(x)$ and $\hat\eta^2(x)$. Two
different gene effect
functions $\alpha(\cdot)$ are implemented. All quantities are
multiplied by 1000} \label{Tb:Comp_SA_NC}
%
\begin{tabular*}{\tablewidth}{@{\extracolsep{\fill}}lcccd{2.2}cd{2.2}@{}}
\hline
&\multicolumn{3}{c}{\textbf{Smooth gene effect}}
&\multicolumn{3}{c@{}}{\textbf{Nonsmooth gene effect}}\\[-4pt]
&\multicolumn{3}{c}{\hrulefill}
&\multicolumn{3}{c@{}}{\hrulefill}
\\
&\multicolumn{1}{c}{$\mathbf{Bias}^{\mathbf{2}}$}
&\multicolumn{1}{c}{$\mathbf{VAR}$}
&\multicolumn{1}{c}{$\mathbf{MISE}$}
&\multicolumn{1}{c}{$\mathbf{Bias}^{\mathbf{2}}$}
&\multicolumn{1}{c}{$\mathbf{VAR}$}
&\multicolumn{1}{c@{}}{$\mathbf{MISE}$}\\
\hline
$\hat{\xi}^2(x)$&0.01&0.14&0.15&16.00&1.47&17.47\\
$\hat{\eta}^2(x)$&0.57&0.24&0.80&0.00&0.22&0.23\\
\hline
\end{tabular*}
\end{table}

In Table \ref{Tb:Comp_SA_NC}, we report the mean integrated squared
bias ($\mathrm{Bias}^2$), the mean integrated variance (VAR),
and the mean integrated squared error (MISE) of $\hat\xi^2(x)$ and
$\hat\eta^2(x)$
with and without the smoothness
assumption on the gene effect $\alpha_g$. From the left panel of
Table \ref{Tb:Comp_SA_NC}, we can see that when the smoothness
assumption is valid, the estimator $\hat\xi^2(x)$ outperforms $\hat
\eta^2(x)$.
The reason is that the mean function $\alpha(X_{\gi})$ depends on the
replication and is not a constant. Therefore, model (\ref{a1}) fails
and $\hat\eta^2(x)$ is biased. One should compare the results with those
on the second row of the right panel where the model is right for
$\hat\eta^2(x)$. In this case, $\hat\eta^2(x)$ performs much better.
Its variance
is about $3/2$ as large as the variance in the case that mean is
generated from a smooth function $\alpha(X_{\gi})$. This is
expected. In the latter case, to eliminate $\alpha_g$, the degree of
freedom reduces from $I=3$ to 2, whereas in the former case,
$\alpha(X_{\gi})$ can be estimated without losing the degree of
freedom, namely the number of replications is still 3. The ratio
$3/2$ is reflected in Table \ref{Tb:Comp_SA_NC}. However, when the smoothness
assumption does not hold, there is serious bias in the estimator
$\hat\xi^2(x)$, even though that $\alpha_g$ is still reasonably sparse.
The bias is an order of magnitude larger than those in the other situations.

%
\begin{table}
\caption{Mean integrated squared bias ($\mathrm{Bias}^2$), mean
integrated variance (VAR), mean integrated squared error (MISE) over
1000 simulations for different variance estimators $\hat\eta^2(x)$ and
$\hat\sigma_O^2(x)$. Seven different correlation
schemes are simulated: $\rho=-0.4$, $\rho=-0.2$, $\rho=0$,
$\rho=0.2$, $\rho=0.4$, $\rho=0.6$ and $\rho=0.8$. All quantities
are multiplied by 1000} \label{Tb:Comp_NC_WC_AGG}
\begin{tabular*}{\tablewidth}{@{\extracolsep{\fill}}lccccccd{2.2}d{2.2}@{}}
\hline
&&\multicolumn{7}{c@{}}{$\bolds\rho$}\\[-4pt]
&&\multicolumn{7}{c@{}}{\hrulefill}\\
&&\multicolumn{1}{c}{$\bolds{-0.4}$}
&\multicolumn{1}{c}{$\bolds{-0.2}$} & \multicolumn{1}{c}{$\bolds{0}$}
&\multicolumn{1}{c}{$\bolds{0.2}$} & \multicolumn{1}{c}{$\bolds{0.4}$}
&\multicolumn{1}{c}{$\bolds{0.6}$} & \multicolumn{1}{c@{}}{$\bolds{0.8}$}\\
\hline
$\mathrm{Bias}^2$&$\hat\eta^2(x)$
&5.93&1.48&0.00&1.48&5.91&13.31&23.67\\[2pt]
&$\hat\sigma_A^2(x)$&0.00&0.00&0.00&0.00&0.00&0.00&0.00\\[2pt]
&$\hat\sigma_O^2(x)$&0.00&0.00&0.00&0.00&0.00&0.00&0.01\\[2pt]
VAR&$\hat\eta
^2(x)$&0.44&0.33&0.24&0.16&0.10&0.05&0.02\\[2pt]
&$\hat\sigma_A^2(x)$&0.27&0.25&0.24&0.22&0.20&0.19&0.20\\[2pt]
&$\hat\sigma_O^2(x)$&0.27&0.25&0.24&0.22&0.20&0.18&0.23\\[2pt]
MISE&$\hat\eta
^2(x)$&6.37&1.81&0.24&1.64&6.01&13.37&23.69\\[2pt]
&$\hat\sigma_A^2(x)$&0.27&0.25&0.24&0.22&0.21&0.19&0.20\\[2pt]
&$\hat\sigma_O^2(x)$&0.27&0.25&0.24&0.22&0.20&0.18&0.24\\
\hline
\end{tabular*}
\end{table}

%
\begin{figure}[b]

\includegraphics{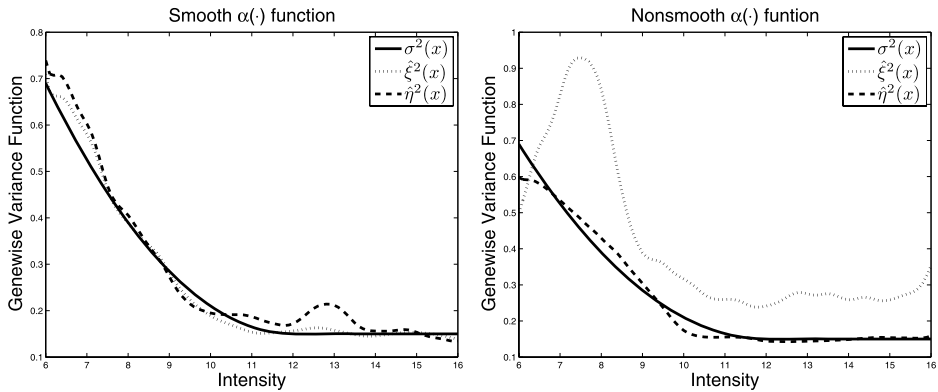}

\caption{Variance estimators $\hat\xi{}^2(x)$ and $\hat\eta^2(x)$
with median performance when different gene effect function
$\alpha(\cdot)$ are implemented. Left panel: smooth
$\alpha(\cdot)$ function. Right panel: nonsmooth $\alpha(\cdot)$
function.} \label{Fig:examAss}
\end{figure}

To see how variance estimators behave, we plot typical estimators
$\hat\xi^2(x)$ and $\hat\eta^2(x)$ with median ISE value among 100
simulations in
Figure \ref{Fig:examAss}. The solid line is the true variance
function while the dotted and dashed lines represent $\hat\xi^2(x)$ and
$\hat\eta^2(x)$, respectively. On the left panel of Figure
\ref{Fig:examAss}, we can see that estimator $\hat\xi^2(x)$ outperforms
the estimator $\hat\eta^2(x)$ when the smoothness assumption is valid. The
region where the biases occur has already been explained above.
However, $\hat\xi^2(x)$ will generate substantial bias when the
nonparametric regression model does not hold, and at the same time,
our nonparametric estimator $\hat\eta^2(x)$ corrects the bias very well.

\subsection{Performance of new estimators}\label{sec33}

In this example, we consider the setting in Section \ref{sec31} that the
smoothness assumption of the gene effect $\alpha_g$ is not valid. For
comparison purpose only, we add an oracle estimator $\hat\sigma_O^2(x)$
in which we assume that $\sigma_1$, $\sigma_2$ and $\rho$ are all known.
We now\vspace*{1pt} evaluate the performance of the estimators $\hat\eta^2(x)$,
$\hat\sigma_A^2(x)$ and $\hat\sigma_O^2(x)$ when the correlation between
within-array replications varies. To be more specific, seven
different correlation settings are considered: $\rho=-0.4, -0.2, 0$,
0.2, 0.4, 0.6, 0.8, with $\rho= 0$ representing across-array
replications. In this case, we increase the number of
simulations to $T=1000$. Again, we report $\mathrm{Bias}^2$, VAR and
MISE of the three estimators for each correlation setting in
Table~\ref{Tb:Comp_NC_WC_AGG}. When $\rho=0$, all the three estimators
give the same bias and variance. This is consistent with our theory.
We can see clearly from the table that, when $\rho\neq0$, the
estimator $\hat\sigma_A^2(x)$ produces much smaller biases than $\hat
\eta^2(x)$. In
fact, when $|\rho|$ as small as 0.2, the bias of $\hat\eta^2(x)$ already
dominates the variance.

%

%
\begin{table}[b]
\caption{Squared bias, variance and MSE of $\hat\rho$,
$\hat\sigma_1$ and $\hat\sigma_2$ in the estimate $\hat\sigma
_A^2(x)$.\break All
quantities are multiplied by $10^6$}
\label{Tb:rho_WC_AGG}
\begin{tabular*}{\tablewidth}{@{\extracolsep{\fill}}lcd{2.2}d{2.2}d{2.2}d{2.2}d{2.2}d{2.2}d{2.2}@{}}
\hline
&&\multicolumn{7}{c@{}}{$\bolds\rho$}\\[-4pt]
&&\multicolumn{7}{c@{}}{\hrulefill}\\
\multicolumn{1}{@{}l}{$\bolds{\hat\sigma_A^2(x)}$}&
&\multicolumn{1}{c}{$\bolds{-0.4}$}&\multicolumn{1}{c}{$\bolds{-0.2}$}
&\multicolumn{1}{c}{$\bolds{0}$}&\multicolumn{1}{c}{$\bolds{0.2}$}
&\multicolumn{1}{c}{$\bolds{0.4}$}&\multicolumn{1}{c}{$\bolds{0.6}$}
&\multicolumn{1}{c@{}}{$\bolds{0.8}$}\\
\hline
$\hat\rho$&$\mathrm{Bias}^2$&0.07&0.04&0.01&0.00&0.00&0.00&3.90\\
&VAR&7.90&16.91&28.65&36.17&35.68&27.21&20.44\\
&MSE&7.97&16.95&28.66&36.17&35.68&27.21&24.35\\
[2pt]
$\hat\sigma_1$&$\mathrm{Bias}^2$&0.24&0.23&0.19&0.14&0.11&0.05&2.47\\
&VAR&11.65&11.52&11.79&12.46&13.64&15.55&18.66\\
&MSE&11.89&11.75&11.99&12.60&13.75&15.59&21.12\\
[2pt]
$\hat\sigma_2$&$\mathrm{Bias}^2$&0.14&0.14&0.12&0.09&0.08&0.05&0.67\\
&VAR&10.34&10.17&10.45&11.12&12.24&13.96&16.16\\
&MSE&10.47&10.31&10.57&11.20&12.32&14.00&16.83\\
\hline
\end{tabular*}
\end{table}

It is worth noticing that the performance of $\hat\sigma_O^2(x)$ and
$\hat\sigma_A^2(x)$ are almost always the same, which indicates that our
algorithm for estimating $\rho$, $\sigma_1$ and $\sigma_2$ is very
accurate. To see this more clearly, the squared bias, variance and
MSE of the estimator $\rho$, $\sigma_1$ and $\sigma_2$ in $\hat
\sigma_A^2(x)$
under the seven correlation settings are reported in Table
\ref{Tb:rho_WC_AGG}. Here, the true value of $\sigma_1$ and
$\sigma_2$ is 0.4217 and 0.1857. For example, when $\rho=0.8$, the
bias of $\hat\rho$ is less than $0.002$ for $\hat\sigma_A^2(x)$,
which is
acceptable because the convergence threshold in the algorithm is set
to be $0.001$.

%
\begin{figure}

\includegraphics{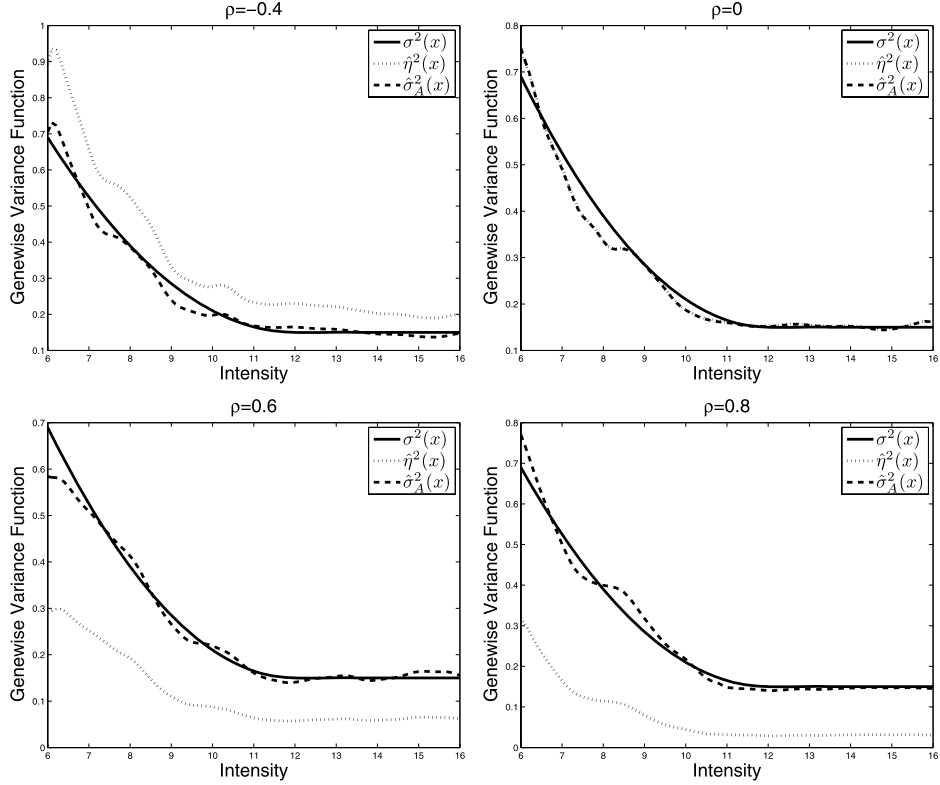}

\caption{Median performance of variance estimators $\hat\eta^2(x)$,
$\hat\sigma^2(x)$ and $\hat\sigma_A^2(x)$ when $\rho=-0.4$, 0,
0.6 and 0.8.}
\label{Fig:Considering_rho}
\end{figure}

In Figure \ref{Fig:Considering_rho}, we render the estimates
$\hat\eta^2(x)$ and $\hat\sigma_A^2(x)$ with the median ISE under
four different
correlation settings: $\rho=-0.4$, $\rho=0$, $\rho=0.6$ and $\rho=
0.8$. We omit the other correlation schemes since they all have
similar performance. The solid lines represent the true variance
function. The dotted lines and dashed lines are for $\hat\eta^2(x)$ and
$\hat\sigma_A^2(x)$, respectively. For the case $\rho=0$, the two estimators
are indistinguishable. When $\rho< 0$, $\hat\eta^2(x)$ overestimates the
genewise variance function, whereas when $\rho> 0$, it
underestimates the genewise variance function.

\section{Application to human total RNA samples using Agilent arrays}\label{sec4}

Our real data example comes from Microarray Quality Control (MAQC)
project [\citet{Patterson}]. The main purpose of the original
paper is on comparison of reproducibility, sensitivity and
specificity of microarray measurements across different platforms
(i.e., one-color and two-color) and testing sites. The MAQC project
use two RNA samples, Stratagene Universal Human Reference total RNA
and Ambion Human Brain Reference total RNA. The two RNA samples have
been assayed on three kinds of arrays: Agilent, CapitalBio and
TeleChem. The data were collected at five sites. Our study focuses
only on the Agilent arrays. At each site, 10 two-color Agilent
microarrays are assayed with 5 of them dye swapped, totaling 30
microarrays.

\subsection{Validation test}\label{Application I}

In the first application, we revisit the validation test as
considered in \citet{FanNiu07}. For the purpose of the validation
tests, we use gProcessedSignal and rProcessedSignal values from
Agilent Feature Extraction software as input. We follow the
preprocessing scheme described in \citet{Patterson} and get 22,144
genes from a total of 41,675 noncontrol genes. Among those, 19 genes
with each having 10 replications are used for validation tests.
Under the null hypothesis of no experimental biases, a reasonable
model is
%
%
\begin{equation} \label{d1}\quad
Y_{\gi} = \alpha_{g} + \varepsilon_{\gi},\qquad \varepsilon_{\gi}\sim
N(0,\sigma^2_g),\qquad
i = 1,\ldots, I, g = 1,\ldots, G.
\end{equation}
We use the notation $G$ to denote the number of genes that have $I$
replications. For our data, $G=19$ and $I=10$. Note that $G$ can be
different from $N$, the total number of different genes. The
validation test statistics in \citet{FanNiu07} include weighted
statistics
\[
T_1=\sum_{g=1}^{G} \Biggl\{
\sum_{i=1}^{I}(Y_{\gi}-\bar{Y}_{g})^2\Big/\sigma_g^2
\Biggr\},\qquad T_2 = \sum_{g=1}^{G} \Biggl\{ {\sum_{i=1}^{I}}
| Y_{\gi}-\bar{Y}_{g} |\Big/\sigma_g \Biggr\},
\]
and unweighted test statistics
\begin{eqnarray*}
T_3&=& \Biggl\{\sum_{g=1}^{G}
\sum_{i=1}^{I}(
Y_{\gi}-\bar{Y}_{g})^2-(I-1)\sum_{g=1}^{G}
\sigma_{g}^2 \Biggr\} \Biggl\{2(I-1)\sum_{g=1}^G \sigma^4_{g} \Biggr\}^{-1/2},\\
T_4&=& \Biggl\{{\sum_{g=1}^{G} \sum_{i=1}^{I}}
|Y_{\gi}-\bar{Y}_{g} |- \lambda_I \sum_{g=1}^{G}
\sigma_{g} \Biggr\} \Bigg/ \Biggl\{ \kappa_I
\Biggl(\sum_{g=1}^{G} \sigma_{g}^2 \Biggr)^{1/2} \Biggr\},
\end{eqnarray*}
where $\lambda_I = \sqrt{2 I(I-1)/\pi}$ and $\kappa_I^2 =
\operatorname{var}({\sum_{i=1}^I} | \varepsilon_{\gi}-\bar{\varepsilon
}_{g} |/
\sigma_g )$.
Under the null hypothesis, the test statistic $T_1$ is $\chi^2$
distributed with degree of freedom
$(I-1)G$ and $T_2, T_3$ and $T_4$ are
all asymptotically normally distributed. As a result, the
corresponding $p$-values can be easily computed.

Here, we apply the same statistics $T_1$, $T_2$, $T_3$ and $T_4$ but
we replace the pooled sample variance estimator by the aggregated
local linear estimator
\[
\hat\sigma^2_g=\sum_{i=1}^I
\hat\sigma_A^2(X_{\gi})\hat{f}(X_{\gi})\Big/\sum_{i=1}^I \hat{f}(X_{\gi}),
\]
where $\hat{f}$ is the estimated density function of $X_{\gi}$. The
difference between the new variance estimator and the simple pooled
variance estimator is that we consider the genewise variance as a
nonparametric function of the intensity level. The latter estimator
may drag small variances of certain arrays to much higher levels by
averaging, resulting in a larger estimated genewise variance and
smaller test statistics or bigger $p$-values.

%
\begin{table}
\tablewidth=280pt
\caption{Comparison of $p$-values for $T_1,\ldots, T_4$ for MAQC
project data considering all 30 arrays together}
\label{Tb: validation_test}
\begin{tabular*}{\tablewidth}{@{\extracolsep{\fill}}lcccc@{}}
\hline
& \multicolumn{4}{c@{}}{$\bolds p$\textbf{-values}} \\[-4pt]
& \multicolumn{4}{c@{}}{\hrulefill} \\
\textbf{Slide name} & $\bolds{T_1}$ & $\bolds{T_2}$
& $\bolds{T_3}$ & $\bolds{T_4}$ \\
\hline
AGL-1-C1 & 1.0000 & 1.0000 & 1.0000 & 1.0000 \\
AGL-1-C2 & 1.0000 & 1.0000 & 1.0000 & 1.0000 \\
AGL-1-C3 & 1.0000 & 1.0000 & 1.0000 & 1.0000 \\
AGL-1-C4 & 1.0000 & 1.0000 & 1.0000 & 1.0000 \\
AGL-1-C5 & 1.0000 & 1.0000 & 1.0000 & 1.0000 \\
AGL-1-D1 & 1.0000 & 1.0000 & 1.0000 & 1.0000 \\
AGL-1-D2 & 1.0000 & 1.0000 & 1.0000 & 1.0000 \\
AGL-1-D3 & 1.0000 & 1.0000 & 1.0000 & 1.0000 \\
AGL-1-D4 & 1.0000 & 1.0000 & 1.0000 & 1.0000 \\
AGL-1-D5 & 1.0000 & 1.0000 & 1.0000 & 1.0000 \\
AGL-2-C1 & 1.0000 & 1.0000 & 1.0000 & 1.0000 \\
AGL-2-C2 & 1.0000 & 1.0000 & 1.0000 & 1.0000 \\
AGL-2-C3 & 1.0000 & 1.0000 & 1.0000 & 1.0000 \\
AGL-2-C4 & 1.0000 & 1.0000 & 1.0000 & 1.0000 \\
AGL-2-C5 & 1.0000 & 1.0000 & 1.0000 & 1.0000 \\
AGL-2-D1 & 1.0000 & 0.9999 & 0.9996 & 0.9999 \\
AGL-2-D2 & 0.8387 & 0.9011 & 0.8953 & 0.9182 \\
AGL-2-D3 & 0.3525 & 0.1824 & 0.3902 & 0.1905 \\
AGL-2-D4 & 1.0000 & 1.0000 & 1.0000 & 1.0000 \\
AGL-2-D5 & 0.8820 & 0.8070 & 0.8848 & 0.7952 \\
AGL-3-C1 & 1.0000 & 1.0000 & 1.0000 & 1.0000 \\
AGL-3-C2 & 1.0000 & 1.0000 & 1.0000 & 1.0000 \\
AGL-3-C3 & 1.0000 & 1.0000 & 1.0000 & 1.0000 \\
AGL-3-C4 & 1.0000 & 1.0000 & 1.0000 & 1.0000 \\
AGL-3-C5 & 1.0000 & 1.0000 & 1.0000 & 1.0000 \\
AGL-3-D1 & 1.0000 & 1.0000 & 1.0000 & 1.0000 \\
AGL-3-D2 & 1.0000 & 1.0000 & 1.0000 & 1.0000 \\
AGL-3-D3 & 1.0000 & 1.0000 & 1.0000 & 1.0000 \\
AGL-3-D4 & 1.0000 & 1.0000 & 1.0000 & 1.0000 \\
AGL-3-D5 & 1.0000 & 1.0000 & 1.0000 & 1.0000 \\
\hline
\end{tabular*}
\end{table}

%
\begin{table}
\tablewidth=280pt
\caption{Comparison of $p$-values for $T_1,\ldots, T_4$ for MAQC
project data considering the arrays with and without dye-swap
separately} \label{Tb:sepTest}
\begin{tabular*}{\tablewidth}{@{\extracolsep{\fill}}lcccc@{}}
\hline
& \multicolumn{4}{c@{}}{$\bolds p$\textbf{-values}} \\[-4pt]
& \multicolumn{4}{c@{}}{\hrulefill} \\
\textbf{Slide name} & $\bolds{T_1}$ & $\bolds{T_2}$
& $\bolds{T_3}$ & $\bolds{T_4}$ \\
\hline
AGL-1-C1 & 1.0000 & 1.0000 & 1.0000 & 1.0000 \\
AGL-1-C2 & 1.0000 & 1.0000 & 1.0000 & 1.0000 \\
AGL-1-C3 & 1.0000 & 1.0000 & 0.9999 & 1.0000 \\
AGL-1-C4 & 1.0000 & 1.0000 & 1.0000 & 1.0000 \\
AGL-1-C5 & 1.0000 & 1.0000 & 0.9999 & 1.0000 \\
AGL-1-D1 & 1.0000 & 1.0000 & 1.0000 & 1.0000 \\
AGL-1-D2 & 1.0000 & 1.0000 & 1.0000 & 1.0000 \\
AGL-1-D3 & 1.0000 & 1.0000 & 1.0000 & 1.0000 \\
AGL-1-D4 & 1.0000 & 1.0000 & 1.0000 & 1.0000 \\
AGL-1-D5 & 1.0000 & 1.0000 & 1.0000 & 1.0000 \\
AGL-2-C1 & 1.0000 & 1.0000 & 0.9943 & 1.0000 \\
AGL-2-C2 & 1.0000 & 1.0000 & 1.0000 & 1.0000 \\
AGL-2-C3 & 1.0000 & 1.0000 & 1.0000 & 1.0000 \\
AGL-2-C4 & 0.0152 & 0.9493 & 0.3931 & 0.9136 \\
AGL-2-C5 & 1.0000 & 1.0000 & 0.8060 & 1.0000 \\
AGL-2-D1 & 0.7806 & 0.8074 & 0.6622 & 0.6584 \\
AGL-2-D2 & 0.2170 & 0.2984 & 0.1651 & 0.2217 \\
AGL-2-D3 & \textbf{0.0002} & \textbf{0.0000}& \textbf{0.0001} &
\textbf{0.0000} \\
AGL-2-D4 & 1.0000 & 1.0000 & 1.0000 & 1.0000 \\
AGL-2-D5 &\textbf{0.1236} &\textbf{0.0662} &\textbf{0.0669}
&\textbf{0.0300} \\
AGL-3-C1 & 1.0000 & 1.0000 & 0.9996 & 1.0000 \\
AGL-3-C2 & 1.0000 & 1.0000 & 0.9988 & 1.0000 \\
AGL-3-C3 & 1.0000 & 1.0000 & 0.9977 & 1.0000 \\
AGL-3-C4 & 1.0000 & 1.0000 & 1.0000 & 1.0000 \\
AGL-3-C5 & 1.0000 & 1.0000 & 0.9999 & 1.0000 \\
AGL-3-D1 & 1.0000 & 1.0000 & 1.0000 & 1.0000 \\
AGL-3-D2 & 1.0000 & 1.0000 & 1.0000 & 1.0000 \\
AGL-3-D3 & 1.0000 & 1.0000 & 1.0000 & 1.0000 \\
AGL-3-D4 & 1.0000 & 1.0000 & 1.0000 & 1.0000 \\
AGL-3-D5 & 1.0000 & 1.0000 & 1.0000 & 1.0000 \\\hline
\end{tabular*}
\end{table}

In the analysis here, we first consider all thirty arrays. The
estimated correlation among replicated genes is $\hat{\rho} = 0.69$.
The $p$-values based on the newly estimated genewise variance are
depicted in Table \ref{Tb: validation_test}. As explained in
\citet{FanNiu07}, $T_4$ is the most stable test among the four.
It turns out that none of the arrays needs further normalization,
which is the same as \citet{FanNiu07}. Furthermore, we separate the
analysis into two groups: the first group using 15 arrays without
dye-swap, which has the estimated correlation $\hat{\rho} = 0.66$,
and the second group using 15 arrays with dye-swap, resulting in an
estimated correlation $\hat{\rho} = 0.34$. The $p$-values are
summarized in Table \ref{Tb:sepTest}. Results show that array
AGL-2-D3 and array AGL-2-D5 need further normalization if 5\%
significance level applies. The difference is due to decreased
estimated $\rho$ for the dye swap arrays and $p$-values are
sensitive to the genewise variance. We also did analysis by
separating data into 6 groups: with and without dye swap, and three
sites of experiments. Due to the small sample size, the six
estimates of $\rho$ range from $0.08$ to $0.74$, and we also find
that array AGL-2-D3 needs further normalization.

\subsection{Gene selection}\label{sec42}

To detect the differentially expressed genes, we follow the filter
instruction and get 19,802 genes out of 41,000 unique
noncontrol genes as in \citet{Patterson}, that is, $I = 1$. The
dye swap result was averaged before doing the one-sample $t$-test.
Thus, at each site, we
have five microarrays.

%
%
%

%
\begin{table}[b]
\tabcolsep=0pt
\caption{Comparison of the number of significantly differentially
expressed genes} \label{Tb:num_significantly_expressed}
\begin{tabular*}{\tablewidth}{@{\extracolsep{4in minus 4in}}lld{5.0}d{5.0}d{5.0}d{5.0}d{5.0}d{5.0}d{5.0}r@{}}
\hline
&&
\multicolumn{2}{c}{$\bolds{p<0.05}$} & \multicolumn{2}{c}{$\bolds{p<0.01}$}
&
\multicolumn{2}{c}{$\bolds{p<0.005}$} & \multicolumn{2}{c@{}}{$\bolds{p<0.001}$}
\\[-4pt]
&&
\multicolumn{2}{c}{\hrulefill} & \multicolumn{2}{c}{\hrulefill}
&
\multicolumn{2}{c}{\hrulefill} &
\multicolumn{2}{c@{}}{\hrulefill}\\
&&\multicolumn{1}{c}{$\bolds t$\textbf{-test}}
&\multicolumn{1}{c}{$\bolds{z}$\textbf{-test}}
& \multicolumn{1}{c}{$\bolds{t}$\textbf{-test}}
& \multicolumn{1}{c}{$\bolds{z}$\textbf{-test}}
& \multicolumn{1}{c}{$\bolds{t}$\textbf{-test}}
& \multicolumn{1}{c}{$\bolds{z}$\textbf{-test}}
& \multicolumn{1}{c}{$\bolds{t}$\textbf{-test}}
&\multicolumn{1}{c@{}}{$\bolds{z}$\textbf{-test}} \\
\hline
Agilent 1 & FC${}>{}$1.5 & 12692 & 12802 & 12464 & 12752 & 12313 & 12722 & 11744 &
12646 \\
& FC${}>{}$2 & 8802 & 8879 & 8654 & 8872 & 8556 & 8869 & 8231 &
8858 \\
& FC${}>{}$4 & 3493 & 3493 & 3431 & 3493 & 3376 & 3493 & 3231 & 3493 \\
Agilent 2 & FC${}>{}$1.5 & 12282 & 12678 & 11217 & 12587 & 10502 & 12536 & 8270 &
12421 \\
& FC${}>{}$2 & 8644 & 8877 & 7908 & 8875 & 7452 & 8861 & 6125 &
8828 \\
& FC${}>{}$4 & 3600 & 3649 & 3188 & 3649 & 2964 & 3649 & 2422 & 3649 \\
Agilent 3 & FC${}>{}$1.5 & 12502 & 12692 & 11994 & 12576 & 11694 & 12519 & 10788 &
12374 \\
& FC${}>{}$2 & 8689 & 8832 & 8344 & 8810 & 8150 & 8800 & 7591 &
8762 \\
& FC${}>{}$4 & 3585 & 3603 & 3378 & 3602 & 3278 & 3602 & 2985 & 3600 \\
\hline
\end{tabular*}
\end{table}

For each site, significant genes are selected based on these 5
dye-swaped average arrays. For all $N = 19\mbox{,}802$ genes, there are no
within-array replications. However, model (\ref{a1}) is still
reasonable, in which $i$ indexes the array. Hence, the
``within-array correlation'' becomes ``between-array correlation''
and is reasonably assumed as $\rho= 0$.

In our nonparametric estimation for the variance function, all the
19,802 genes are used to estimate the variance function, which gives
us enough reason to believe that the estimator $\hat\sigma_A^2(x)$ is
close to
the inherent true variance function $\sigma^2(x)$.

We applied both the $t$-test and $z$-test to each gene to see if the
logarithm of the expression ratio is zero, using the five arrays
collected at each location. The number of differentially expressed
genes detected by using the two different tests under three Fold
Changes (FC) and four significant levels are given in Table
\ref{Tb:num_significantly_expressed}. Large numbers of genes are
identified as differentially expressed, which is expected when
comparing a brain sample and a tissue pool sample
[\citet{Patterson}]. We can see clearly that the $z$-test
associated with our new variance estimator $\hat\sigma_A^2(x)$ leads
to more
differentially expressed genes. For example, at site~1, using
$\alpha= 0.001$, among the fold changes at least 2, $t$-test picks
8231 genes whereas $z$-test selects 8875 genes. This gives an
empirical power increase of $(8875-8231)/19\mbox{,}802 \approx3.25\%$ in
the group with observed fold change at least~2.

To verify the accuracy of our variance estimation in the $z$-test,
we compare the empirical power increase with the expected
theoretical power increase. The expected theoretical power increase is computed
as
%
%
\begin{equation} \label{d2}
\mathrm{ave}\{\mathrm{P}_z (\mu_g/\sigma_g) -
\mathrm{P}_{t_{n-1}}(\mu_g/\sigma_g)\},
\end{equation}
taking the average of power increases across all $\mu_g \not= 0$.
However, in the absence of the availability, we replace $\mu_g$ by
its estimate, which is the sample average of $n=5$ observed
log-expression ratios. Table \ref{Tb:power_diff} depicts the
results at three different sites, in which the columns ``Theo''
refer to the expected theoretical power increase defined by
(\ref{d2}), with $\mu_g$ replaced by $\bar{Y}_g$ and $\sigma_g$
replaced by its estimate from the genewise variance function, and
the columns ``Emp'' refer to the empirical power increase.

%
\begin{table}
\caption{Comparison of expected theoretical and empirical power difference
(in percentage)}\label{Tb:power_diff}
\begin{tabular*}{\tablewidth}{@{\extracolsep{\fill}}lccd{2.2}d{2.2}d{2.2}d{2.2}d{2.2}c@{}}
\hline
& \multicolumn{2}{c}{$\bolds{\alpha=0.05}$}
& \multicolumn{2}{c}{$\bolds{\alpha=0.01}$}
& \multicolumn{2}{c}{$\bolds{\alpha=0.005}$}
& \multicolumn{2}{c@{}}{$\bolds{\alpha=0.001}$}\\[-4pt]
&
\multicolumn{2}{c}{\hrulefill} & \multicolumn{2}{c}{\hrulefill}
&
\multicolumn{2}{c}{\hrulefill} &
\multicolumn{2}{c@{}}{\hrulefill}\\
& \multicolumn{1}{c}{\textbf{Theo}} & \multicolumn{1}{c}{\textbf{Emp}}
& \multicolumn{1}{c}{\textbf{Theo}} & \multicolumn{1}{c}{\textbf{Emp}}
& \multicolumn{1}{c}{\textbf{Theo}} & \multicolumn{1}{c}{\textbf{Emp}}
& \multicolumn{1}{c}{\textbf{Theo}} & \multicolumn{1}{c@{}}{\textbf{Emp}}\\
\hline
Agilent 1 & 2.52 & 0.61 & 6.08 & 3.75 & 8.06 & 5.59 & 13.66 & 11.74 \\
Agilent 2 & 4.03 & 7.56 & 10.11 & 17.61 & 13.61 & 22.86 & 23.75 & 37.63\\
Agilent 3 & 3.02 & 2.56 & 7.14 & 7.39 & 9.42 & 10.19 & 15.94 & 18.18
\\[2pt]
Average & 3.19 & 3.58 & 7.78 & 9.58 & 10.36 & 12.88 & 17.79 & 22.51 \\
\hline
%
%
\end{tabular*}
\end{table}

There are two things worth noticing here. First, for expected theoretical
power increase, we use the sample mean $\bar{Y}_g = \mu_g +
\bar\epsilon_g$ instead of the real gene effect~$\mu_g$, which is
not observable, so it inevitably involves the error term
$\bar\epsilon_g$. Second, the power functions $\mathrm{P}_z(\mu)$ and
$\mathrm{P}_t(\mu)$ depend sensitively on $\mu$ and the tails of the
assumed distribution. Despite of these, the expected theoretical and
empirical power increases are in the same bulk and
the averages are very close.
This provides good evidence that our genewise variance estimation
has an acceptable accuracy.

We also apply SIMEX and permutation SIMEX methods in \citet{Carroll08}
to the MAQC data, to
illustrate its utility. As mentioned in the \hyperref[sec1]{Introduction}, their model
is not really intended for
the analysis of two-color microarray data. Should we only use the
information on log-ratios ($Y$), the
model is very hard to interpret. In addition, one might question why
the information on $X$
(observed intensity levels) is not used at all. Nevertheless, we apply
the SIMEX methods of \citet{Carroll08} to only the log-ratios $Y$
in the two-color data and
produce similar tables to the
Tables \ref{Tb:num_significantly_expressed} and \ref{Tb:power_diff}.

From the results, we have the following understandings. First, all the
numbers for $z$-test in
Tables \ref{Tb:num_significantly_expressed-SIMEX} and \ref
{Tb:num_significantly_expressed-perSIMEX}
at all significance levels are approximately the same. In fact, the
$p$-values are very small so that numbers at
different significance levels are about the same. That indicates that
both SIMEX and permutation SIMEX
method are tending to estimate genewise variance very small, making the
test statistics large for
all the time. On the other hand, our method estimates the genewise
variance moderately so that the
numbers are not exactly the same for different significance levels.
Second, in the implementation,
we found that the SIMEX and permutation SIMEX is computationally
expensive (more than one hour)
while our method only takes a few minutes. Third, from Tables \ref{Tb:
power_diff-SIMEX} and
\ref{Tb: power_diff-perSIMEX} we can see that the expected theoretical
power increase and the
empirical ones are reasonably close, which are in lines with our method.

%
\begin{table}
\tabcolsep=0pt
\caption{Comparison of the number of significantly differentially
expressed genes using SIMEX method}
\label{Tb:num_significantly_expressed-SIMEX}
\begin{tabular*}{\tablewidth}{@{\extracolsep{4in minus 4in}}lld{5.0}d{5.0}d{5.0}d{5.0}d{5.0}d{5.0}d{5.0}r@{}}
\hline
&&
\multicolumn{2}{c}{$\bolds{p<0.05}$} & \multicolumn{2}{c}{$\bolds{p<0.01}$}
&
\multicolumn{2}{c}{$\bolds{p<0.005}$} & \multicolumn{2}{c@{}}{$\bolds{p<0.001}$}
\\[-4pt]
&&
\multicolumn{2}{c}{\hrulefill} & \multicolumn{2}{c}{\hrulefill}
&
\multicolumn{2}{c}{\hrulefill} &
\multicolumn{2}{c@{}}{\hrulefill}\\
&&\multicolumn{1}{c}{$\bolds t$\textbf{-test}}
&\multicolumn{1}{c}{$\bolds{z}$\textbf{-test}}
& \multicolumn{1}{c}{$\bolds{t}$\textbf{-test}}
& \multicolumn{1}{c}{$\bolds{z}$\textbf{-test}}
& \multicolumn{1}{c}{$\bolds{t}$\textbf{-test}}
& \multicolumn{1}{c}{$\bolds{z}$\textbf{-test}}
& \multicolumn{1}{c}{$\bolds{t}$\textbf{-test}}
&\multicolumn{1}{c@{}}{$\bolds{z}$\textbf{-test}} \\
\hline
Agilent 1 & FC${}>{}$1.5 & 12692 & 12820 & 12464 & 12820 & 12313 & 12820 & 11744 &
12820 \\
& FC${}>{}$2 & 8802 & 8879 & 8654 & 8879 & 8556 & 8879 & 8231 &
8879 \\
& FC${}>{}$4 & 3493 & 3493 & 3431 & 3493 & 3376 & 3493 & 3231 & 3493 \\
Agilent 2 & FC${}>{}$1.5 & 12282 & 12721 & 11217 & 12721 & 10502 & 12721 & 8270 &
12721 \\
& FC${}>{}$2 & 8644 & 8878 & 7908 & 8878 & 7452 & 8878 & 6125 &
8878 \\
& FC${}>{}$4 & 3600 & 3649 & 3188 & 3649 & 2964 & 3649 & 2422 & 3649 \\
Agilent 3 & FC${}>{}$1.5 & 12502 & 12760 & 11994 & 12760 & 11694 & 12760 & 10788 &
12760 \\
& FC${}>{}$2 & 8689 & 8836 & 8344 & 8836 & 8150 & 8836 & 7591 &
8836 \\
& FC${}>{}$4 & 3585 & 3603 & 3378 & 3603 & 3278 & 3603 & 2985 & 3603 \\
\hline
\end{tabular*}
\end{table}

%
\begin{table}
\tabcolsep=0pt
\caption{Comparison of the number of significantly differentially
expressed genes using permutation SIMEX}
\label{Tb:num_significantly_expressed-perSIMEX}
\begin{tabular*}{\tablewidth}{@{\extracolsep{4in minus 4in}}lld{5.0}d{5.0}d{5.0}d{5.0}d{5.0}d{5.0}d{5.0}r@{}}
\hline
&&
\multicolumn{2}{c}{$\bolds{p<0.05}$} & \multicolumn{2}{c}{$\bolds{p<0.01}$}
&
\multicolumn{2}{c}{$\bolds{p<0.005}$} & \multicolumn{2}{c@{}}{$\bolds{p<0.001}$}
\\[-4pt]
&&
\multicolumn{2}{c}{\hrulefill} & \multicolumn{2}{c}{\hrulefill}
&
\multicolumn{2}{c}{\hrulefill} &
\multicolumn{2}{c@{}}{\hrulefill}\\
&&\multicolumn{1}{c}{$\bolds t$\textbf{-test}}
&\multicolumn{1}{c}{$\bolds{z}$\textbf{-test}}
& \multicolumn{1}{c}{$\bolds{t}$\textbf{-test}}
& \multicolumn{1}{c}{$\bolds{z}$\textbf{-test}}
& \multicolumn{1}{c}{$\bolds{t}$\textbf{-test}}
& \multicolumn{1}{c}{$\bolds{z}$\textbf{-test}}
& \multicolumn{1}{c}{$\bolds{t}$\textbf{-test}}
&\multicolumn{1}{c@{}}{$\bolds{z}$\textbf{-test}} \\
\hline
Agilent 1 & FC${}>{}$1.5 & 12692 & 12820 & 12464 & 12820 & 12313 & 12820 & 11744 &
12820 \\
& FC${}>{}$2 & 8802 & 8879 & 8654 & 8879 & 8556 & 8879 & 8231 &
8879 \\
& FC${}>{}$4 & 3493 & 3493 & 3431 & 3493 & 3376 & 3493 & 3231 & 3493 \\
Agilent 2 & FC${}>{}$1.5 & 12282 & 12721 & 11217 & 12721 & 10502 & 12721 & 8270 &
12721 \\
& FC${}>{}$2 & 8644 & 8878 & 7908 & 8878 & 7452 & 8878 & 6125 &
8878 \\
& FC${}>{}$4 & 3600 & 3649 & 3188 & 3649 & 2964 & 3649 & 2422 & 3649 \\
Agilent 3 & FC${}>{}$1.5 & 12502 & 12760 & 11994 & 12760 & 11694 & 12760 & 10788 &
12760 \\
& FC${}>{}$2 & 8689 & 8836 & 8344 & 8836 & 8150 & 8836 & 7591 &
8836 \\
& FC${}>{}$4 & 3585 & 3603 & 3378 & 3603 & 3278 & 3603 & 2985 & 3603 \\
\hline
\end{tabular*}
\end{table}

%
\begin{table}
\caption{Comparison of expected theoretical and empirical power
difference using SIMEX method
(in percentage)}\label{Tb: power_diff-SIMEX}
\begin{tabular*}{\tablewidth}{@{\extracolsep{\fill}}lccd{2.2}d{2.2}d{2.2}d{2.2}d{2.2}c@{}}
\hline
& \multicolumn{2}{c}{$\bolds{\alpha=0.05}$}
& \multicolumn{2}{c}{$\bolds{\alpha=0.01}$}
& \multicolumn{2}{c}{$\bolds{\alpha=0.005}$}
& \multicolumn{2}{c@{}}{$\bolds{\alpha=0.001}$}\\[-4pt]
&
\multicolumn{2}{c}{\hrulefill} & \multicolumn{2}{c}{\hrulefill}
&
\multicolumn{2}{c}{\hrulefill} &
\multicolumn{2}{c@{}}{\hrulefill}\\
& \multicolumn{1}{c}{\textbf{Theo}} & \multicolumn{1}{c}{\textbf{Emp}}
& \multicolumn{1}{c}{\textbf{Theo}} & \multicolumn{1}{c}{\textbf{Emp}}
& \multicolumn{1}{c}{\textbf{Theo}} & \multicolumn{1}{c}{\textbf{Emp}}
& \multicolumn{1}{c}{\textbf{Theo}} & \multicolumn{1}{c@{}}{\textbf{Emp}}\\
\hline
Agilent 1 & 2.43 & 2.06 & 7.17 & 5.42 & 10.30 & 7.34 & 20.71 & 13.44 \\
Agilent 2 & 7.16 & 3.41 & 19.20 & 12.06 & 26.17 & 16.90 & 43.46 & 30.42
\\
Agilent 3 & 4.18 & 2.88 & 11.71 & 7.38 & 16.45 & 9.89 & 31.38 & 17.57
\\ [2pt]
Average & 4.59 & 2.78 & 12.69 & 8.29 & 17.64 & 11.38 & 31.85 & 20.48 \\
\hline
\end{tabular*}
\end{table}

%
\begin{table}
\caption{Comparison of expected theoretical and empirical power
difference using permutation\break SIMEX method
(in percentage)}\label{Tb: power_diff-perSIMEX}
\begin{tabular*}{\tablewidth}{@{\extracolsep{\fill}}lccd{2.2}d{2.2}d{2.2}d{2.2}d{2.2}c@{}}
\hline
& \multicolumn{2}{c}{$\bolds{\alpha=0.05}$}
& \multicolumn{2}{c}{$\bolds{\alpha=0.01}$}
& \multicolumn{2}{c}{$\bolds{\alpha=0.005}$}
& \multicolumn{2}{c@{}}{$\bolds{\alpha=0.001}$}\\[-4pt]
&
\multicolumn{2}{c}{\hrulefill} & \multicolumn{2}{c}{\hrulefill}
&
\multicolumn{2}{c}{\hrulefill} &
\multicolumn{2}{c@{}}{\hrulefill}\\
& \multicolumn{1}{c}{\textbf{Theo}} & \multicolumn{1}{c}{\textbf{Emp}}
& \multicolumn{1}{c}{\textbf{Theo}} & \multicolumn{1}{c}{\textbf{Emp}}
& \multicolumn{1}{c}{\textbf{Theo}} & \multicolumn{1}{c}{\textbf{Emp}}
& \multicolumn{1}{c}{\textbf{Theo}} & \multicolumn{1}{c@{}}{\textbf{Emp}}\\
\hline
Agilent 1 & 1.89 & 2.86 & 5.66 & 6.43 & 8.19 & 8.59 & 16.75 & 15.07 \\
Agilent 2 & 4.84 & 7.37 & 13.44 & 17.22 & 18.97 & 22.50 & 36.90 & 37.26
\\
Agilent 3 & 2.89 & 4.91 & 8.34 & 10.13 & 11.87 & 13.11 & 23.44 & 21.31
\\ [2pt]
Average & 3.20 & 5.05 & 9.15 & 11.26 & 13.01 & 14.74 & 25.70 & 24.55 \\
\hline
\end{tabular*}
\end{table}

\section{Summary}\label{sec5}

The estimation of genewise variance function is motivated by the
downstream analysis of microarray data such as validation test and
selecting statistically differentially expressed genes. The
methodology proposed here is novel by using across-array and within-array
replications. It does not require any specific assumptions on
$\alpha_g$ such as sparsity or smoothness, and hence reduces the
bias of the conventional nonparametric estimators. Although the
number of nuisance parameters is proportional to the sample size, we
can estimate the main interest (variance function) consistently. By
increasing the degree of freedom largely, both the validation tests
and $z$-test using our variance estimators are more powerful in
identifying arrays that need to be normalized further and more
capable of selecting differentially expressed genes.

Our proposed methodology has a wide range of applications. In
addition to the microarray data analysis with within-array
replications, it can be also applied to the case without within-array
replications, as long as the model (\ref{a1}) is reasonable.
Our two-way nonparametric model is a natural extension of the
Neyman--Scott problem. Therefore, it is applicable to all the
problems where the Neyman--Scott problem is applicable.

There are possible extensions. For example, the SIMEX idea can be
applied on our model in order to take into account the measurement error.
We can also make adaptations to our methods when we have a prior
correlation structure among replications other than the identical
correlation assumption.


\begin{appendix}\label{app}
\section*{Appendix}

The following regularity conditions are imposed for the technical
proofs:
\begin{enumerate}
\item[1.] The regression function $\sigma^2(x)$ has
a bounded and continuous second derivative.

\item[2.]The kernel function $K$ is a bounded symmetric density
function with a
bounded support.

\item[3.]$h \rightarrow0, Nh \rightarrow\infty. $

\item[4.]$\mathrm{E}[\sigma^8(X)]$ exists and the marginal density
$f_X(\cdot)$ is continuous.
\end{enumerate}

We need the following conditional variance--covariance matrix of the
random vector ${\mathbf{Z}_g}$ in our asymptotic study.
\begin{lemma}\label{marvar}
Let $\bolds{\Omega}$ be the variance--covariance matrix of ${\mathbf{Z}_g}$
conditioning on all data $\mathbf{X}$. Then, respectively, the
diagonal and
off-diagonal elements are
%
%
%
\begin{eqnarray}\label{app1}
\Omega_{ii} &=& 2\sigma^4(X_{\gi}) + \frac{2}{(I-1)^2(I-2)^2}\sum_{k
\neq l} \sigma^2(X_{gk}) \sigma^2(X_{gl}) \nonumber\\[-8pt]\\[-8pt]
&&{} + \frac{4(I-3)}{(I-1)(I-2)^2}\sigma^2(X_{\gi}) \sum_{j \neq
i} \sigma^2(X_{gj}),\qquad i = 1,\ldots, I,\nonumber
\\
\label{app2}
\Omega_{ij} &=& \frac{4}{(I-1)^2} \sigma^2(X_{\gi}) \sigma^2(X_{gj})\nonumber\\
&&{}
+ \frac{2}{(I-1)^2(I-2)^2}\mathop{\sum_{k \neq l}}_{k,l \neq i,j} \sigma
^2(X_{gk}) \sigma^2(X_{gl})\nonumber\\[-8pt]\\[-8pt]
&&{} - \frac{4}{(I-1)^2(I-2)}\sum_{k \neq i,j} \sigma^2(X_{gk}) \bigl(\sigma
^2(X_{\gi})+ \sigma^2(X_{gj}) \bigr), \nonumber\\
\eqntext{i,j = 1,\ldots, I, i \neq j.}
\end{eqnarray}
\end{lemma}
\begin{pf}
Let $\mathbf{A}$ be the variance--covariance matrix of ${\mathbf{r}_g}$
conditioning on all data~$\mathbf{X}$. By direct computation, the diagonal
elements are given by
%
%
\begin{eqnarray}\label{app3}
\hspace*{22pt}A_{ii} &=& \operatorname{var}[(Y_{\gi}-\bar{Y}_g)^2|{\mathbf
{X}}]\nonumber\\
\hspace*{22pt}&=& \frac{2(I-1)^4}{I^4}\sigma^4(X_{\gi})+\frac{4(I-1)^2}{I^4} \sum_{k
\neq i}\sigma^2 (X_{\gi}) \sigma^2(X_{gk})
+\frac{2}{I^4}
\sum_{k \neq i} \sigma^4(X_{gk})\\
\hspace*{22pt}&&{}+\frac{4}{I^4} \sum_{l,k
\neq i, l < k} \sigma^2(X_{gl})\sigma^2(X_{gk}),\qquad i = 1,\ldots, I,\nonumber
\end{eqnarray}
and the off-diagonal elements are given by
%
%
\begin{eqnarray}\label{app4}
A_{ij} &=& \operatorname{cov}\{[(Y_{\gi}-\bar{Y}_g)^2,(Y_{gj}-\bar
{Y}_g)^2]|{\mathbf{X}}\}
\nonumber\\
&=&\frac{2(I-1)^2}{I^4}[\sigma^4(X_{\gi})+\sigma^4(X_{gj})] +\frac
{4(I-1)^2}{I^4}\sigma^2(X_{\gi})\sigma^2(X_{gj})\nonumber\\[-8pt]\\[-8pt]
&&{} -\frac{4(I-1)}{I^4}\sum_{k \neq
i,j}\sigma^2(X_{gk})\bigl(\sigma^2(X_{\gi}) + \sigma^2(X_{gj})\bigr)\nonumber\\
&&{} +\frac{4}{I^4}\sum_{k,l \neq
i,j;l<k}\sigma^2(X_{gl})\sigma^2(X_{gk})+\frac{2}{I^4}\sum_{k \neq
i,j} \sigma^4(X_{gk}).\nonumber
\end{eqnarray}
Using $\bolds{\Omega}= \mathbf{B}\mathbf{A}\mathbf{B}^T$, we can
obtain the result by direct
computation.
\end{pf}

The proofs for Theorems \ref{T1} and \ref{T3} follow a similar idea. Since
Theorem \ref{T1} does not involve a lot of coefficients, we will show the
proof of Theorem \ref{T1} and explain the difference in the proof of
Theorem \ref{T3}.
\begin{pf*}{Proof of Theorem \ref{T1}}
First of all, the bias of $\eta_i^2(x)$ comes from the local linear
approximation. Since $\{(X_{\gi}, Z_{\gi})\}_{g=1}^N$ is an i.i.d.
sequence, by (\ref{b2}) and the result in \citet{Fan96}, it follows
that
\[
\mathrm{E}\{ \eta_i^2(x) | \mathbf{X}\} = \sigma^2(x) + b(x) + o_P(h^2).
\]
Similarly, the asymptotic variance of $\eta_i^2(x)$ also follows
from \citet{Fan96}.

We now prove the off-diagonal elements in matrix $\operatorname
{var}[\bolds{\eta}|\mathbf{X}]$
%
%
\begin{equation}\label{app5}
\operatorname{cov}[(\hat\eta_i^2(x),\hat\eta_j^2(x))|{\mathbf
{X}}] =V_2+o_P(1/N).
\end{equation}
%
Recalling that $\hat\eta^2_i(x) = \sum_{g=1}^N
W_{N,i}((X_{\gi}-x)/h)Z_{\gi}$, we have
\[
\operatorname{cov}[(\hat\eta_i^2(x),\hat\eta_j^2(x))|{\mathbf
{X}}] = \sum_{g=1}^N
W_{N,i} \biggl( \frac{X_{\gi}-x}{h} \biggr) W_{N,j} \biggl( \frac{X_{gj}-x}{h} \biggr)
\operatorname{cov}
[(Z_{\gi},Z_{gj})| \mathbf{X}].
\]
The equality follows by the fact that $\operatorname
{cov}[(Z_{\gi},Z_{g'j})| \mathbf{X}] =
0$ when $g \neq g'$. Recall $\Omega_{ij} = \operatorname
{cov}[(Z_{\gi}, Z_{gj})|\mathbf{X}
]$ and define $R_{N,g} = N\cdot W_{N,j}((X_{gj}-x)/h) \Omega_{ij}$. Thus,
%
%
\begin{equation}\label{app6}
N\cdot\operatorname{cov}[(\hat\eta_i^2(x),\hat\eta
_j^2(x))|{\mathbf{X}}] = \sum
_{g=1}^N W_{N,i} \biggl( \frac{X_{\gi}-x}{h} \biggr) R_{N,g}.
\end{equation}
The right-hand side of (\ref{app6}) can be seen as local linear
smoother of the synthetic data $\{(X_{\gi}, R_{N,g})\}_{g=1}^N$.
Although $R_{N,g}$ involves $N$ at the first glance, its conditional
expectation $\mathrm{E}[R_{N,g}|X_{\gi}=x]$ and conditional
variance\break
$\operatorname{var}
[R_{N,g}|X_{\gi}=x]$ do not grow with $N$. Since $\{(X_{\gi}, R_{N,g})\}
_{g=1}^N$ is an i.i.d. sequence, by the results in \citet{Fan96},
we obtain
\[
N\cdot\operatorname{cov}[(\hat\eta_i^2(x),\hat\eta
_j^2(x))|{\mathbf{X}}] = \mathrm{E}
[R_{N,g}|X_{\gi}=x] + o_P(1).
\]
To calculate $\mathrm{E}[R_{N,g}|X_{\gi}=x]$, we apply the approximation
$W_{N,i}(u) = K(u)(1+o_P(1))/(Nhf_X(x))$ in the example of Fan and
Gijbels [(\citeyear{Fan96}), page 64] and have the following arguments
\begin{eqnarray*}
&&\mathrm{E}[R_{N,g}|X_{\gi}=x] \\
&&\qquad= \mathrm{E}\biggl[ N \cdot\frac{1}{Nhf_X(x)}
hK_h(X_{gj}-x) \Omega_{ij} |X_{\gi}=x \biggr]\bigl(1+o_P(1)\bigr) \\
&&\qquad=(f_X(x))^{-1} \int K(u) \Omega_{ij}|_{X_{\gi}=x}(x+hu,\mathbf{s})
f_X(x+hu)\,du\,d\mathbf{s}+ o_P(1)\\
&&\qquad= N V_2 + o_P(1),
\end{eqnarray*}
where $\mathbf{s}$ represents all the integrating variables
corresponding to
$X_{g1},\ldots, X_{gI}$ except $X_{\gi}$ and $X_{gj}$. That justifies
(\ref{app5}).

To prove the multivariate asymptotic normality
%
%
\begin{equation}\label{app7}
\bolds\Sigma^{-1/2}\bigl(\bolds{\eta}-\bigl(\sigma^2(x)+b(x)+o_P(h^2)\bigr)
\mathbf{e}\bigr)
\stackrel{D}\longrightarrow N(0, \mathbf{I}_I),
\end{equation}
we employ Cram\'er--Wold device: for any unit vector $\mathbf{a}=
(a_1,\ldots, a_I)^T$ in $\mathbb{R}^I$,
\[
F^* \stackrel{\triangle}{=}\{\mathbf{a}^T \bolds\Sigma\mathbf{a}\}
^{-1/2} \Biggl\{ \sum_{i=1}^I a_i \sum_{g=1}^N W_{N,i} \biggl( \frac
{X_{\gi}-x}{h} \biggr)\bigl(Z_{\gi}-\sigma^2(Z_{\gi})\bigr) \Biggr\} \stackrel
{D}\longrightarrow N(0,1).
\]
%
Denote by $Q_{g,i} = W_{N,i}((X_{\gi}-x)/h)(Z_{\gi} - \sigma^2(X_{\gi}))$
and $\widetilde{Q}_g = \sum_{i=1}^I a_i Q_{g,i} $. Note that the
sequence $\{
\widetilde{Q}_g\}_{g=1}^N$ is i.i.d. distributed. To show the asymptotic
normality of $F^*$, it is sufficient to check Lyapunov's condition:
\[
\lim_{N\rightarrow\infty} \frac{\sum_{g=1}^N \mathrm
{E}[|\widetilde{Q}_g|^4 | \mathbf{X}
] }{(\sum_{g=1}^N \mathrm{E}[ |\widetilde{Q}_g|^2 | \mathbf{X}]
)^2} = 0.
\]
To facilitate the presentation, we first note that sequences $\{
Q_{g,i}\}_{g=1}^N$ are i.i.d. and satisfy Lyapunov's condition for each
fixed $i$. Denote $\delta^2_{N,i} =\break \sum_{g=1}^N \mathrm
{E}[|Q_{g,i}|^2|\mathbf{X}
]$. And recall that $\delta^2_{N,i} = \operatorname{var}[\hat\eta
^2_i(x)| \mathbf{X}] =
O_P((Nh)^{-1})$. Let $c^*$ be a generic constant which may vary from
one line to another. We have the following approximation:
\begin{eqnarray*}
\sum_{g=1}^N \mathrm{E}[|Q_{g,i}|^4 | \mathbf{X}] &=& c^* N^{-3}
\mathrm{E}\bigl\{
K_h^4(X_{\gi}-x)\bigl[\bigl(Z_{\gi}-\sigma^2(X_{\gi})\bigr)^4|\mathbf{X}\bigr] \bigr\}\bigl(1+o_P(1)\bigr)
\\
&=& O_P((Nh)^{-3}).
\end{eqnarray*}
Therefore, $\sum_{g=1}^N \mathrm{E}[|Q_{g,i}|^4|\mathbf{X}] =
o(\delta^4_{N,i})$.
By the marginal Lyapunov conditions, we have the following inequality:
\[
\sum_{g=1}^N\mathrm{E}[\widetilde{Q}_g^4 |\mathbf{X}] \leq c^* \sum
_{i=1}^I \sum_{g=1}^N
\mathrm{E}[|Q_{g,i}|^4|\mathbf{X}]
= c^*I \cdot o_P((Nh)^{-2}) = o_P((Nh)^{-2}).
\]
For the denominator, we have the following arguments:
\begin{eqnarray*}
\sum_{g=1}^N \mathrm{E}[|\widetilde{Q}_g|^2|\mathbf{X}] 
&=&\sum_{i}a^2_i \sum_{g=1}^N \mathrm{E}[Q^2_{g,i}|\mathbf{X}] +
\sum_{i \neq j}
a_i a_j \sum_{g=1}^N \mathrm{E}[Q_{g,i} Q_{g,j}| \mathbf{X}]\\
&=& \sum_{i}a_i^2 \operatorname{var}[\hat\eta^2_i(x)|\mathbf{X}] +
\sum_{i \neq j} a_i
a_j \operatorname{cov}[(\hat\eta^2_i(x),\hat\eta^2_j(x))|\mathbf
{X}]\\
&\stackrel{*}{=}& O_P((Nh)^{-1}) + O_P(N^{-1}) \\
&=& O_P((Nh)^{-1}).
\end{eqnarray*}
Note that the second to last equality holds by the asymptotic
conditional variance--covariance matrix $\bolds\Sigma$. Therefore Lyapunov's
condition is justified. That completes the proof.
\end{pf*}
%
%
\begin{pf*}{Proof of Theorem \ref{P1}}
First of all, for each given $g$,
\[
\mathrm{E}s^2_{B,g} = I \operatorname{var}(\bar{Y}_{gj}) = \sigma_2 +
\rho(I-1)\sigma_1^2.
\]
Note that by (\ref{b6}), we have
\begin{eqnarray*}
\mathrm{E} (Y_{\mathit{gij}} - \bar{Y}_{gj})^2 & = & I^{-2} [I(I-1) \sigma_2 + \rho
(I-1)(I-2) \sigma_1^2 - 2(I-1)^2\rho\sigma_1^2]\\
& = & I^{-1}(I-1) ( \sigma_2 - \rho\sigma_1^2).
\end{eqnarray*}
Thus, for all $g$, we have
\[
\mathrm{E}s^2_{W,g} = \sigma_2 - \rho\sigma_1^2.
\]
Since $\{s_{B, g}^2\}$ and $\{s_{W, g}^2\}$ are i.i.d. sequences
across the $N$ genes, by the central limit theorem, we have
\begin{eqnarray*}
\frac{1}{N}\sum_{g=1}^N s^2_{B,g} &=& \sigma_2 + \rho(I-1)\sigma
_1^2 + O_P(N^{-1/2}),\\
\frac{1}{N} \sum_{g=1}^N s^2_{W,g} &=& \sigma_2 - \rho\sigma_1^2 +
O_P(N^{-1/2}).
\end{eqnarray*}
Therefore,
\begin{eqnarray*}
\hat\rho_0 & = & \frac{\sigma_2 + \rho(I-1) \sigma_1^2 - \sigma_2
+ \rho\sigma_1^2 + O_P(N^{-1/2})}{\sigma_2 + \rho(I-1) \sigma_1^2
+ (I-1)(\sigma_2 - \rho\sigma_1^2) + O_P(N^{-1/2})}\\
& = & \rho\sigma_1^2 / \sigma_2 + O_P(N^{-1/2}).
\end{eqnarray*}
\upqed\end{pf*}
\begin{pf*}{Proof of Theorem \ref{T3}}
Note that
\begin{eqnarray*}
\operatorname{var}[\hat\eta_A^2(x)| \mathbf{X}] &=& \sum
_{g=1}^N\sum_{i=1}^I W^2_N
\biggl(\frac{X_{\gi}-x}{h} \biggr)\operatorname{var}[Z_{\gi}|\mathbf{X}]\\
&&{} + \sum_{g=1}^N \sum_{i \neq j}^I W_N \biggl(\frac{X_{\gi}-x}{h} \biggr)W_N
\biggl(\frac{X_{gj}-x}{h} \biggr)\operatorname{cov}[(Z_{\gi},Z_{gj})|\mathbf{X}].
\end{eqnarray*}

Following similar steps in the proof of Theorem \ref{T1}, one can
verify $\operatorname{var}[\hat\eta_A^2(x)|\break \mathbf{X}]=V'_1/I +
(1-1/I)V_2' +
o_P((Nh)^{-1})$, where the coefficients $C_2,\ldots, C_4, D_0$,\break $\ldots,
D_4$ are as follows:
\begin{eqnarray*}
C_2 &=& \frac{4(1+\rho^2)\sigma_2 + [4\rho(I-2)+4\rho
^2(2I-3)]\sigma^2_1}{I-1},\\
C_3 &=& -\frac{8\rho^2(I-3)\sigma^3_1 + 8(\rho^2+\rho)\sigma_1
\sigma_2}{I-1},\\
C_4 &=& \frac{2}{(I-1)(I-2)} \{ (1+\rho^2)\sigma_2^2+2(\rho^2+\rho
)(I-3)\sigma_1^2\sigma_2\\
&&\hspace*{133.3pt}{}+(I-3)(I-4)\rho^2\sigma_1^4 \},
\\
D_0 &=& 2 \biggl( \rho^2- \frac{4\rho}{I-1} + \frac{2(1+\rho^2)}{(I-1)^2}
\biggr),\\
D_1 &=& \frac{8}{(I-1)^2} \{ (2I-4)\rho- (I^2-4I+5)\rho^2 \}\sigma
_1,\\
D_2 &=& \frac{4}{(I-1)^2(I-2)} \bigl\{ (I-3)^2\rho^2+\bigl((I-2)^2+1\bigr)\rho
-2(I-2) \bigr\}\sigma_2 \\
&&{} + \frac{4(I-3)}{(I-1)^2(I-2)} \bigl\{ \bigl(3(I-2)(I-3)+2\bigr)\rho^2-2(I-2)\rho
\bigr\}\sigma_1^2,\\
D_3 &=& -\frac{8(I-3)^2}{(I-1)^2(I-2)} \{ (\rho^2+\rho)\sigma
_1\sigma_2 +(I-4)\rho^2\sigma_1^3 \},\\
D_4 &=& \frac{4}{(I-1)^2(I-2)^2}\\
&&{}\times \biggl\{ (1+\rho^2) \pmatrix{I-2 \cr2}
\sigma^2_2 \\
&&\hspace*{17pt}{}+6(\rho^2+\rho) \pmatrix{I-2 \cr3}\sigma_1^2\sigma_2
+12\rho^2 \pmatrix{I-2 \cr4} \sigma_1^4 \biggr\}.
\end{eqnarray*}
\upqed\end{pf*}
\end{appendix}

\section*{Acknowledgments}
The authors thank the Editor, the Associate Editor and two referees,
whose comments have greatly improved the scope and presentation of
the paper.

\printaddresses

\end{document}